\documentclass[12pt,reqno]{amsart}
\usepackage{hyperref,amssymb,mathrsfs,eqnarray,xcolor,bbold,amsxtra}
\theoremstyle{plain}
\newtheorem{theorem}{Theorem}[section]
\newtheorem{proposition}[theorem]{Proposition}
\newtheorem{corollary}[theorem]{Corollary}
\newtheorem{lemma}[theorem]{Lemma}

\theoremstyle{definition}

\newtheorem{example}[theorem]{Example}

\numberwithin{equation}{section}
\DeclareMathOperator{\Span}{span}
\DeclareMathOperator{\ran}{ran}
\DeclareMathOperator{\rad}{Rad}
\DeclareMathOperator{\dist}{dist}
\DeclareMathOperator{\ins}{ins}
\renewcommand{\le}{\leqslant}
\renewcommand{\ge}{\geqslant}

\begin{document}
\title[On operators with many spectral components]{On operators with many spectral components}
\author[J. Bra\v{c}i\v{c}]{Janko Bra\v{c}i\v{c}}
\address{Faculty of Natural Sciences and Engineering, University of Ljubljana, 
A\v{s}ker\v{c}eva c. 12, SI-1000 Ljubljana, Slovenia}
\email{janko.bracic@ntf.uni-lj.si}
\author[M. Kandi\'c]{Marko Kandi\'c}
\address{Faculty of Mathematics and Physics, University of Ljubljana,
Jadranska 19, 1000 Ljubljana, Slovenia \ \ \ and \ \ \
Institute of Mathematics, Physics, and Mechanics, Jadranska 19, 1000 Ljubljana, Slovenia}
\email{marko.kandic@fmf.uni-lj.si}
\begin{abstract}
For a complex Banach algebra $\mathscr A$ with unit and $n\in\mathbb N$, let $\mathcal C_n(\mathscr A)$ be the 
set of all elements whose spectrum has at least $n$ connected components, and let 
$d_\sigma(\mathscr A)=\sup\{n;\;\mathcal C_n(\mathscr A)\neq\varnothing\}$ be the spectral depth of 
$\mathscr A$. We characterize nonemptiness of $\mathcal C_n(\mathscr A)$ in terms of orthogonal idempotents, 
compute $d_\sigma(\mathscr A)$ for finite-dimensional algebras, and show that every class 
$\mathcal C_n(\mathscr A)$ is open. For operators on an infinite-dimensional complex Banach space $X$ these 
classes are dense in $\mathscr B(X)$, which implies that $\mathcal C_\infty(\mathscr B(X))$ is a dense 
$G_\delta$-subset of $\mathscr B(X)$. By contrast, we exhibit a commutative Banach algebra whose spectral 
classes are all nonempty and none of which is dense. For self-adjoint operators on a Hilbert space, we obtain 
sharper perturbation results, with positive semi-definite perturbations of prescribed rank producing spectra with 
exactly $n$ components. We show that a quasinilpotent weighted forward shift $S_w$ on $X$ admits a nilpotent rank-one  operator $F$ such that the spectrum of $S_w+\alpha F$ has infinitely many components for every nonzero 
$\alpha$.
\end{abstract}
\maketitle
%
\section{Introduction} \label{sec:introduction}
The spectrum of an element of a complex Banach algebra is a nonempty compact subset of the complex plane.
Much of spectral theory concerns how a spectrum's shape reflects the properties of the element. Among the 
crudest, and yet most robust, features of a compact set is the number of its connected components. The purpose 
of this paper is to study this number systematically. We organize the elements of a Banach algebra into classes 
according to how many components their spectra have, and we investigate the resulting hierarchy.

Let $\mathscr A$ be a complex Banach algebra with unit $e$. For each $n\in\mathbb N$ we call
$$\mathcal C_n(\mathscr A)=\{a\in\mathscr A;\;\sigma(a)\text{ has at least }n\text{ components}\}$$
the $n$-th \emph{spectral class} of $\mathscr A$, and we set 
$\mathcal C_\infty(\mathscr A)=\bigcap_{n\in\mathbb N}\mathcal C_n(\mathscr A)$. 
These classes form a decreasing chain beginning with $\mathcal C_1(\mathscr A)=\mathscr A$, and the point 
at which the chain becomes empty is measured by the \emph{spectral depth}
$$d_\sigma(\mathscr A)=\sup\{n\in\mathbb N;\;\mathcal C_n(\mathscr A)\neq\varnothing\}\in
\mathbb N\cup\{\infty\}.$$
Several questions arise at once. Which values can $d_\sigma(\mathscr A)$ take, and how is it computed for 
concrete algebras? How are the classes $\mathcal C_n(\mathscr A)$ situated inside $\mathscr A$ topologically? 
The present paper is devoted to these questions and to the perturbation problems they lead to.

The starting point is that the whole hierarchy is governed by idempotents. In 
Proposition~\ref{prop:exact-spectral-components-idempotents} we show that $\mathcal C_n(\mathscr A)$ is 
nonempty precisely when the unit of $\mathscr A$ decomposes into $n$ nonzero pairwise orthogonal idempotents. 
It follows, in particular, that possessing an element whose spectrum has \emph{at least} $n$ components is the 
same as possessing one with \emph{exactly} $n$ components. The corresponding statement at the level 
$n=\infty$ is more delicate, since a decreasing chain of nonempty sets may have empty intersection. 
Theorem~\ref{thm:Cinfty-orthogonal-idempotents} shows that an infinite family of pairwise orthogonal 
idempotents produces, for every sequence of weights in a suitable weighted $\ell^1$-space, an element with 
prescribed spectrum, and in particular an element of $\mathcal C_\infty(\mathscr A)$. These results allow us to 
compute the spectral depth of finite-dimensional algebras from their Wedderburn decomposition 
(Proposition~\ref{prop:depth-finite-dim-algebra}) and to show that $d_\sigma(\mathscr B(X))=\infty$ for every 
infinite-dimensional complex Banach space $X$ (Corollary~\ref{cor:spectral-depth-BX}).

The topological structure of the classes is the subject of Section~\ref{sec:topological-properties-spectral-classes}. 
Theorem~\ref{thm:Ck-open} states that $\mathcal C_n(\mathscr A)$ is open in every complex Banach algebra with 
unit. The proof rests on a continuity property of the Riesz functional calculus, which we prove in 
Lemma~\ref{lem:continuity-functional-calculus}. Whether the classes are also large is a different matter, and 
here the answer depends on the algebra. For $\mathscr A=C(\Omega)$ the classes are governed by the topology of 
$\Omega$ (Propositions \ref{prop:Ck-spectral-depth} and \ref{prop:Cn-dense-totally-disconnected} and
Example \ref{ex:all-spectral-classes-nonempty-none-dense}). In contrast, for operators on 
an infinite-dimensional Banach space, the classes are as large as possible. 
Building on a perturbation result of Herrero and Salinas \cite{HS72}, we prove in Theorem~\ref{thm:Cn-density} that 
every $T\in\mathscr B(X)$ can be moved into $\mathcal C_n(\mathscr B(X))$ by a perturbation of rank at most $n-1$ 
and arbitrarily small norm. Together with openness, this yields, via the Baire category theorem, that 
$\mathcal C_\infty(\mathscr B(X))$ is a dense $G_\delta$-subset of $\mathscr B(X)$ 
(Theorem~\ref{thm:Cinfty-residual}).

The same section also treats quotients. If $\mathscr I$ is a proper closed two-sided ideal of $\mathscr A$, the 
preimages under the quotient map of the classes $\mathcal C_n(\mathscr A/\mathscr I)$ are open, invariant under 
perturbation by elements of $\mathscr I$, and dense exactly when the corresponding classes are dense in the 
quotient (Proposition~\ref{prop:quotient-Ck-basic-properties}). When the quotient is finite-dimensional, its 
Wedderburn decomposition determines everything 
(Proposition~\ref{prop:finite-dimensional-quotient-depth}). Taking $\mathscr I=\mathscr K(X)$, so that the 
quotient is the Calkin algebra, gives statements about the essential spectrum. The finite-dimensionality 
hypothesis cannot be dropped, as Example~\ref{ex:nonempty-essential-spectral-class-not-dense} shows by means of 
a space whose Calkin algebra is isomorphic to $C(\Omega)$ for a disjoint union $\Omega$ of $n$ disks 
\cite{Mot24}.

Sections~\ref{sec:perturbations-Hilbert-space} and \ref{sec:quasinilpotent-rank-one-perturbation} are devoted to 
sharper perturbation results, where the general theory can be improved by exploiting additional structure. On a 
Hilbert space, a self-adjoint operator whose spectrum has $k$ components can be moved into 
$\mathcal C_n(\mathscr B(H))$, for any $n>k$, by a \emph{positive semi-definite} perturbation of rank exactly 
$n-k$ and arbitrarily small norm, and the resulting spectrum has \emph{exactly} $n$ components 
(Theorem~\ref{thm:positive-exact-n}). Iterating the construction produces a positive trace-class perturbation 
of arbitrarily small trace norm after which the spectrum acquires infinitely many components, described 
explicitly (Theorem~\ref{thm:positive-schatten-infinite-components}). Since trace-class operators are compact, 
the essential spectrum is left untouched. For normal operators, the corresponding statement requires compact 
rather than finite-rank perturbations and is obtained from the Weyl--von Neumann--Berg theorem 
(Theorem~\ref{thm:normal-compact-perturbation-Cinfty}). The perturbation cannot in general be taken normal. 
Proposition~\ref{prop:unitary-normal-perturbation-trivial} shows that a unitary operator without eigenvalues 
admits no nonzero compact normal perturbation preserving normality.

Finally, in Section~\ref{sec:quasinilpotent-rank-one-perturbation} we show how far a single rank-one 
perturbation can move an operator within the hierarchy. Let $S_w$ be a weighted forward shift with respect to a 
basic sequence of $X$, with summable nonincreasing positive weights. Such an operator is compact, and it is 
quasinilpotent exactly when the geometric means of the weights tend to zero 
(Lemma~\ref{lem:weighted-compact-quasinilpotent-shift}). Hence, its spectrum $\{0\}$ is as small and as 
connected as possible. Nevertheless, for a large class of rank-one operators $F$ the perturbed operator 
$S_w+\alpha F$ lies in $\mathcal C_\infty(\mathscr B(X))$ for every nonzero $\alpha$ 
(Theorem~\ref{thm:Cinfty-gamma-sequence-criterion}). The proof identifies the nonzero spectrum of 
$S_w+\alpha F$ with the solution set of an equation $G(z)=1/\alpha$ for an entire function $G$ determined by 
the perturbation, and then appeals to the Great Picard theorem.

\section{Preliminaries} \label{sec:preliminaries}
In this section, we collect notation and a few standard results that we need throughout the paper.

\subsection*{Connected components of a topological space}
Let $\Omega$ be an arbitrary topological space. A subset $\Delta$ of $\Omega$ is said to be \emph{clopen} if it 
is both open and closed. If $\Delta$ can be written as a union of two nonempty disjoint relatively open subsets, then 
it is \emph{disconnected}. Otherwise, we call $\Delta$ \emph{connected}. A \emph{connected component} 
or simply a \emph{component} of $\Delta$ is a maximal connected subset of $\Delta$. 

An easy induction and the fact that in a compact Hausdorff space, the component containing a point is the 
intersection of all clopen subsets containing the point (see e.g. \cite[Theorem~6.1.23]{Eng89}) yield the following 
lemma, whose proof we omit.

\begin{lemma}
\label{lem:finite-clopen-spectral-partition}
Let $n\in \mathbb N$ and let $\Omega$ be a compact Hausdorff space with at least $n$ components. Then there exist 
nonempty pairwise disjoint clopen subsets $\Delta_1,\ldots,\Delta_n$ whose union is $\Omega$.\qed
\end{lemma}

\subsection*{Unital complex Banach algebras}
Throughout this paper, let $\mathscr A$ be a complex Banach algebra with unit $e$. As usual, the {\em spectrum} 
of $a\in \mathscr A$ is $\sigma(a)=\{\lambda\in\mathbb C;\; \lambda e-a\;\text{is not invertible in}\;\mathscr A\}$. 
It is well-known that $\sigma(a)$ is a nonempty compact subset of $\mathbb C$. Hence, its complement 
$\rho(a)=\mathbb C\setminus\sigma(a)$, the {\em resolvent set} of $a$, is an unbounded open set. 
By $r(a)=\max\{|\lambda|;\; \lambda\in \sigma(a)\}$ we denote the 
\emph{spectral radius} of $a$. By the Gelfand--Beurling formula, $r(a)=\lim\limits_{n\to\infty}\|a^n\|^{\frac 1n}$. 
An element $a$ is \emph{quasinilpotent} if $r(a)=0$, that is, $\sigma(a)=\{0\}$.

For any finite collection $\Gamma$ of positively oriented Jordan curves contained in $\mathbb{C}$, the \emph{inside} 
of $\Gamma$ is the open set $\ins(\Gamma)$ consisting of those points of $\mathbb{C}$ around which 
$\Gamma$ has winding number one. Let $a\in \mathscr A$ and let $f$ be a holomorphic function on an open 
neighborhood $U$ of $\sigma(a)$. We say that $\Gamma$ encloses $\sigma(a)$ if 
$\Gamma\subseteq U\setminus\sigma(a)$ and $\sigma(a)\subseteq\ins(\Gamma)$. In this case, we can 
associate with $f$ the element
\begin{equation}\label{eq:holomorphic-calculus}
f(a) = \frac{1}{2\pi i} \int_\Gamma f(\lambda)(\lambda e-a)^{-1}\,d\lambda
\end{equation}
in $\mathscr A$, which does not depend on the choice of $\Gamma$ but only on the values of $f$ on $U$.
The mapping $f\mapsto f(a)$ is an algebra homomorphism, called the \emph{Riesz functional calculus}. It follows,
in particular, that elements $f(a)$ and $g(a)$ in $\mathscr A$ commute whenever $f$ and $g$ are
holomorphic functions on an open neighborhood of $\sigma(a)$. The reader is referred to \cite[Section~I.7]{BD73}
and \cite[Section~VII.4]{Con90} for more details about the functional calculus of a single element in a Banach algebra.

Suppose that $\sigma$ is a clopen subset of the spectrum $\sigma(a)$ of $a\in\mathscr A$. Since $\sigma$ 
and $\sigma(a)\setminus\sigma$ are disjoint compact subsets of $\mathbb C$, they admit disjoint open 
neighbourhoods. Hence, the characteristic function $\chi_\sigma$ extends to a locally constant, and therefore 
holomorphic, function on an open neighbourhood of $\sigma(a)$. The element $e_\sigma=\chi_\sigma(a)$ is the 
\emph{spectral idempotent} or simply, the idempotent associated with $\sigma$. Idempotents $e_1,\ldots,e_n$
are \emph{pairwise orthogonal} if $e_ie_j=0$ for $i\ne j$. For the proof of the following proposition, 
see \cite[Proposition~I.7.9]{BD73}.
 
\begin{proposition}\label{prop:orthogonal-idempotent-decomposition}
Let $a$ be an element of a complex Banach algebra $\mathscr A$ with unit $e$. If $\sigma_1,\ldots,\sigma_n$ are 
pairwise disjoint nonempty clopen subsets of $\sigma(a)$ satisfying $\sigma(a)=\sigma_1\cup\cdots\cup\sigma_n$, 
then the spectral idempotents $e_{\sigma_1},\ldots,e_{\sigma_n}$ are nonzero, pairwise orthogonal, and satisfy
$e=e_{\sigma_1}+\cdots+e_{\sigma_n}$.
\end{proposition}

\subsection*{Complex Banach spaces and operators}
Let $X$ be a complex Banach space with dimension $\dim(X)\geq 2$. We denote its topological dual by $X^*$
and the pairing between $X$ and $X^*$ is given by $\langle x,\xi\rangle=\xi(x)$ for all $x\in X$ and $\xi\in X^*$. 
The Banach algebra of all operators on $X$, i.e., bounded linear transformations, is denoted by $\mathscr B(X)$;
the identity operator is denoted by $I$. The closed two-sided ideal of all compact operators in $\mathscr B(X)$ 
is denoted by $\mathscr K(X)$. For $0\ne x\in X$ and $0\ne \xi\in X^*$, the rank-one operator 
$x\otimes \xi\colon X\to X$ is defined by $(x\otimes \xi)(y)=\langle y,\xi\rangle x$ for all $y\in X$. If $x=0$ or $\xi=0$, 
then $x\otimes \xi$ is the zero operator. For each $m\in\mathbb N$, let 
$\mathscr F_m(X)=\{ x_1\otimes \xi_1+\cdots+x_m\otimes \xi_m;\; x_1,\ldots, x_m\in X,\; \xi_1,\ldots, \xi_m\in X^*\}$
be the set of all operators whose rank is at most $m$, and let 
$\mathscr F(X)=\bigcup\limits_{m\in\mathbb N}\mathscr F_m(X)$ be the set of all finite-rank operators. The set 
$\mathscr F(X)$ is an ideal in $\mathscr B(X)$ (nonclosed if $X$ is infinite-dimensional) which is contained in the ideal
of compact operators. Whenever $X$ is infinite-dimensional and $K\in \mathscr K(X)$, every nonzero point of 
$\sigma(K)$ is an isolated eigenvalue of finite multiplicity, and $0$ is the only possible accumulation point of 
$\sigma(K)$. Moreover, since eigenvectors for nonzero eigenvalues lie in the range and the range of a finite-rank 
operator is finite-dimensional, the spectrum of each finite-rank operator is finite. For instance,
$\sigma(x\otimes \xi)=\{0,\langle x,\xi\rangle\}$ (recall that $\dim(X)\geq 2$).

If $T\in \mathscr B(X)$ and $\sigma\subseteq\sigma(T)$ is a clopen subset, then the spectral idempotent 
$P_\sigma=\chi_\sigma(T)$ is usually called a \emph{spectral projection}. Its range $X_\sigma=P_\sigma X$
is the \emph{spectral subspace} of $T$ corresponding to $\sigma$. It is an invariant subspace for $T$.
Assume that $\sigma_1, \ldots, \sigma_k$ are disjoint clopen subsets of $\sigma(T)$ such that 
$\sigma(T)=\sigma_1\cup \cdots\cup \sigma_k$. Then $X=X_{\sigma_1}\oplus \cdots\oplus X_{\sigma_k}$, 
where the sum is algebraic and topological, that is, every vector $x\in X$ has a unique decomposition 
$x=x_1+\cdots+x_k$, where $x_j\in X_{\sigma_j}$ for each $j=1,\ldots, k$, and for each $j=1,\ldots, k$ 
the mapping $x\mapsto x_j$ is continuous. The operator $T$ decomposes as $T=T_1\oplus\cdots\oplus T_k$
with $T_j\in\mathscr B(X_{\sigma_j})$ and $\sigma(T_j)=\sigma_j$ for $j=1,\ldots,k$.

Lastly, we discuss the essential spectrum $\sigma_e(T)$ of an operator $T\in \mathscr B(X)$. 
First, we recall the definition of a \emph{Fredholm operator}. Let the kernel and the range of $T$ be denoted by 
$\ker(T)$ and $\ran(T)$, respectively. Then $T$ is said to be Fredholm if $\ker(T)$ and $X/\ran(T)$ are 
finite-dimensional.

The quotient Banach algebra $\mathcal Q(X)=\mathscr B(X)/\mathscr K(X)$ is called the \emph{Calkin algebra} 
of $X$. By $\pi\colon \mathscr B(X) \to \mathcal Q(X)$ we denote the associated quotient projection. 
The \emph{essential spectrum} $\sigma_e(T)$ of $T\in \mathscr B(X)$ is defined as the spectrum of $\pi(T)$ in 
the Calkin algebra. Hence, $\lambda\in\sigma_e(T)$ if and only if $\pi(\lambda I-T)$ is not invertible in 
$\mathcal Q(X)$. By Atkinson's theorem \cite{Atk51}, a complex number $\lambda$
belongs to $\sigma_e(T)$ if and only if the operator $\lambda I-T$ is not Fredholm.

\section{Spectral classes and spectral depth} \label{sec:spectral-classes-spectral-depth}
For each $n\in \mathbb N$ the $n$-th \emph{spectral class} $\mathcal C_n(\mathscr A)$ is defined by 
$$\mathcal C_n(\mathscr A)= \{a\in\mathscr A;\;\sigma(a)\text{ has at least } n \text{ components}\}.$$
The definition gives the chain $\mathscr A=\mathcal C_1(\mathscr A)\supseteq\mathcal C_2(\mathscr A)\supseteq
\mathcal C_3(\mathscr A)\supseteq\cdots$, where the first equality holds because the spectrum of every element 
is nonempty. The class $\mathcal C_2(\mathscr A)$ consists precisely of those elements whose spectrum is
disconnected. The following proposition provides a characterization of $\mathcal C_n(\mathscr A)$. 

\begin{proposition}
\label{prop:exact-spectral-components-idempotents}
For each $n\in\mathbb N$ the following statements are equivalent:
\begin{enumerate}
\item there exists $a\in\mathscr A$ such that $\sigma(a)$ has exactly $n$ components;
\item $\mathcal C_n(\mathscr A)\neq \varnothing$;
\item there exist nonzero pairwise orthogonal idempotents $e_1,\ldots,e_n\in\mathscr A$ such that
$e=e_1+\cdots+e_n$;
\item the algebra $\mathscr A$ contains a subalgebra which contains $e$ and is isomorphic to $\mathbb C^n$.
\end{enumerate}
\end{proposition}

\begin{proof}
(i)$\Rightarrow$(ii) is clear.

(ii)$\Rightarrow$(iii) Since, by the assumption, there exists $a\in\mathscr A$ such that $\sigma(a)$ has at least $n$
components, by Lemma~\ref{lem:finite-clopen-spectral-partition} there exist pairwise disjoint nonempty clopen
subsets $\sigma_1,\ldots,\sigma_n$ of $\sigma(a)$ such that $\sigma(a)=\sigma_1\cup\cdots\cup\sigma_n$.
The existence of required pairwise orthogonal idempotents follows from
Proposition~\ref{prop:orthogonal-idempotent-decomposition}. 

(iii)$\Rightarrow$(iv) Define a mapping $\Phi\colon \mathbb C^n\to \mathscr A$ by 
$\Phi(\alpha_1,\ldots,\alpha_n)=\sum_{k=1}^n\alpha_k e_k$. Clearly, $\Phi$ is linear and bounded.
Note that $\Phi(1,\ldots,1)=e$. The multiplication in $\mathbb C^n$ is pointwise. Hence, the pairwise
orthogonality of $e_1,\ldots,e_n$ gives 
\begin{align*}
\Phi\big((\alpha_1,\ldots,\alpha_n)(\beta_1,\ldots,\beta_n)\big)&=\sum_{k=1}^n\alpha_k\beta_k e_k
=\big(\sum_{k=1}^n\alpha_k e_k\big)\big(\sum_{k=1}^n\beta_k e_k\big)\\
&= \Phi(\alpha_1,\ldots,\alpha_n)\Phi(\beta_1,\ldots,\beta_n).
\end{align*}
Thus, $\Phi$ is a unital homomorphism. To prove the injectivity, suppose that $\Phi(\alpha_1,\ldots,\alpha_n)=0$, 
that is, $\alpha_1 e_1+\cdots+\alpha_n e_n=0$. Multiplying this equality by $e_j$, and using pairwise orthogonality, 
gives $\alpha_j e_j=0$; since $e_j\ne 0$, we get $\alpha_j=0$. Hence, $\Phi$ is injective.
We conclude that $\Phi$ is an isomorphism from $\mathbb C^n$ onto a subalgebra of $\mathscr A$ that contains $e$.

(iv)$\Rightarrow$(iii) 
Let $\Phi\colon \mathbb C^n\to \mathscr A$ be an isomorphism onto a subalgebra $\mathscr D\subseteq 
\mathscr A$ with $e\in\mathscr D$. Then $e$ is the unit of $\mathscr D$ and therefore $\Phi(1,\ldots,1)=e$.
A direct calculation shows that $\Phi$ maps the standard basis vectors of $\mathbb C^n$ to the required pairwise 
orthogonal idempotents $e_1,\ldots,e_n\in\mathscr A$, which gives $\Phi(1,\ldots,1)=e_1+\cdots+e_n$.
We conclude that $e=e_1+\cdots+e_n$.

(iii)$\Rightarrow$(i) Let $e_1,\ldots,e_n\in\mathscr A$ be nonzero orthogonal idempotents in $\mathscr A$ that 
sum up to $e$. Choose distinct complex numbers $\lambda_1,\ldots,\lambda_n$, and set
$a=\lambda_1e_1+\cdots+\lambda_ne_n$. We claim that $\sigma(a)=\{\lambda_1,\ldots,\lambda_n\}$.
To this end, pick $\lambda\notin\{\lambda_1,\ldots,\lambda_n\}$, and define
$$b = (\lambda_1-\lambda)^{-1}e_1+\cdots+(\lambda_n-\lambda)^{-1}e_n.$$
Using the pairwise orthogonality of the idempotents $e_1,\ldots, e_n$ and the identity $e=e_1+\cdots+e_n$, 
we obtain $(a-\lambda e)b=b(a-\lambda e)=e$, which shows that $\lambda\in\rho(a)$.
On the other hand, fix $1\le j\le n$. Then 
$(a-\lambda_j e)e_j=
((\lambda_1-\lambda_j) e_1+\cdots+(\lambda_j-\lambda_j)e_j+\cdots+(\lambda_n-\lambda_j) e_n)e_j=0$. 
Since $e_j\neq 0$, we see that $a-\lambda_j e$ is not invertible. Thus, $\lambda_j\in\sigma(a)$, proving 
$\sigma(a)=\{\lambda_1,\ldots,\lambda_n\}$. Since $\lambda_1,\ldots,\lambda_n$ are distinct numbers, 
$\sigma(a)$ has exactly $n$ components.
\end{proof}

The equivalence between (i) and (ii) in Proposition~\ref{prop:exact-spectral-components-idempotents} shows that 
for complex Banach algebras with unit there is no difference between possessing elements with ``at least $n$'' and 
``exactly $n$'' spectral components.

\begin{corollary}
Every element in $\mathscr A$ has a connected spectrum if and only if $\mathscr A$ has no idempotents other 
than $0$ and $e$.
\end{corollary}

\begin{proof}
This is the equivalence between negations of items (ii) and (iii) for $n=2$, upon noting that a nontrivial idempotent
$p$ yields the decomposition $e=p+(e-p)$ into nonzero orthogonal idempotents, and conversely.
\end{proof}

\begin{corollary}
\label{cor:Ck-chain-strict}
If for some $n\in \mathbb N$ we have $\mathcal C_n(\mathscr A)\neq \varnothing$, then 
$\mathcal C_{n+1}(\mathscr A) \subsetneq \mathcal C_n(\mathscr A)$.
\end{corollary}

\begin{proof}
Suppose that $\mathcal C_n(\mathscr A)\neq \varnothing$. Then, by 
Proposition \ref{prop:exact-spectral-components-idempotents}, there exist nonzero pairwise orthogonal idempotents 
$e_1,\ldots,e_n\in\mathscr A$ such that $e=e_1+\cdots+e_n$. Now using the part ``(iii)$\Rightarrow$(i)'' in 
the proof of Proposition \ref{prop:exact-spectral-components-idempotents}, one can construct an element 
$a\in\mathscr A$ whose spectrum consists of exactly $n$ components. Hence, 
$a\in \mathcal C_n(\mathscr A)\setminus\mathcal C_{n+1}(\mathscr A)$.
\end{proof}

We now turn to the limiting spectral class corresponding to spectra with infinitely many components. Define
$$\mathcal C_\infty(\mathscr A)=\bigcap_{n=1}^{\infty}\mathcal C_n(\mathscr A).$$
Hence, $\mathcal C_\infty(\mathscr A)$ consists precisely of those elements in $\mathscr A$ 
whose spectra have infinitely many components. The \emph{spectral depth} of $\mathscr A$ is
$$d_\sigma(\mathscr A)=\sup\{n\in\mathbb N;\;\mathcal C_n(\mathscr A)\neq\varnothing\}\in\mathbb N
\cup\{\infty\}.$$
If $d_\sigma(\mathscr A)$ is a positive integer $n$, then Corollary~\ref{cor:Ck-chain-strict} gives
$$\mathscr A=\mathcal C_1(\mathscr A) \supsetneq\cdots\supsetneq \mathcal C_{n}(\mathscr A)
\supsetneq \mathcal C_{n+1}(\mathscr A)=\varnothing$$
and Proposition~\ref{prop:exact-spectral-components-idempotents} immediately implies $\dim(\mathscr A)\geq n$.

The spectral classes have been introduced for an arbitrary complex Banach algebra with a unit, 
but no generality is lost in the operator setting. Indeed, let $\mathscr A$ be a complex Banach 
algebra with unit $e$ and let $L\colon \mathscr A\to \mathscr B(\mathscr A)$ be the 
\emph{left regular representation}, that is, $L_a(x)=ax$ for all $a,x\in\mathscr A$. Then $L$ is 
an algebra homomorphism, $\|L_a\|\leq \|a\|$, and $L$ is injective because $L_a(e)=a$. Moreover, 
$L$ preserves spectra, that is, $\sigma(L_a)=\sigma(a)$ for every $a\in\mathscr A$.
The inclusion $\sigma(L_a)\subseteq\sigma(a)$ is clear, since $L_{a^{-1}}$ is the inverse of $L_a$ 
whenever $a$ is invertible. Conversely, assume that $L_a$ is invertible in 
$\mathscr B(\mathscr A)$ and put $b=L_a^{-1}(e)$. Then $ab=L_a(b)=e$, and therefore 
$L_aL_b=L_{ab}=L_e=I$. Since $L_a$ is invertible, we get $L_b=L_a^{-1}$, hence 
$L_{ba}=L_bL_a=I=L_e$, and injectivity of $L$ gives $ba=e$. Thus, $a$ is invertible in 
$\mathscr A$. Consequently, $a\in\mathcal C_n(\mathscr A)$ if and only if 
$L_a\in\mathcal C_n(\mathscr B(\mathscr A))$ for every $n\in\mathbb N$; intersecting over $n$, the same 
equivalence holds for $\mathcal C_\infty(\mathscr A)$.
Let us emphasize that injectivity alone does not suffice; for a general injective homomorphism $\pi$ one has 
only $\sigma(\pi(a))\subseteq \sigma(a)$.

We nevertheless retain the Banach algebra formulation. The main reason is that the inequality 
$d_\sigma(\mathscr A)\leq d_\sigma(\mathscr B(\mathscr A))$ carries no information when 
$\dim(\mathscr A)=\infty$, because $d_\sigma(\mathscr B(X))=\infty$ for every infinite-dimensional 
$X$ as we shall see in Corollary~\ref{cor:spectral-depth-BX}, whereas the spectral depth of $\mathscr A$ itself 
may well be finite; for instance, $d_\sigma\big(C[0,1]\big)=1$. 

The condition $d_\sigma(\mathscr A)=\infty$ asserts that every class $\mathcal C_n(\mathscr A)$ is nonempty, 
whereas $\mathcal C_\infty(\mathscr A)\neq\varnothing$ asserts that their intersection is nonempty. The former 
does not formally imply the latter, since a decreasing chain of nonempty sets may have empty intersection; 
whether the implication nevertheless holds is studied in \cite{BK26}. For this reason we now exhibit elements 
with infinitely many spectral components directly. The construction rests on an infinite family of pairwise 
orthogonal idempotents, and it prescribes the resulting spectrum completely.

\begin{theorem}
\label{thm:Cinfty-orthogonal-idempotents}
Let $(p_j)_{j\in\mathbb N}$ be a sequence of nonzero pairwise orthogonal idempotents in $\mathscr A$.
For every $j\in\mathbb N$, put $w_j=\|p_j\|$ and let 
$$\ell^1(w)=\Big\{(\lambda_j)_{j\in\mathbb N}\subseteq\mathbb C;\;
\sum_{j=1}^\infty|\lambda_j|\,w_j<\infty\Big\}.$$
Then for every $(\lambda_j)_{j\in\mathbb N}\in\ell^1(w)$ the series $\sum_{j=1}^\infty\lambda_jp_j$ converges 
absolutely in $\mathscr A$ and its sum $a$ satisfies
$$\sigma(a)=\{0\}\cup\{\lambda_j;\; j\in\mathbb N\}.$$
If, in addition, the numbers $\lambda_j$ are pairwise distinct and nonzero, then $a\in\mathcal C_\infty(\mathscr A)$.
\end{theorem}

\begin{proof}
Since $0\ne p_j=p_j^2$, we have $0<\|p_j\|\leq\|p_j\|^2$, and therefore $w_j=\|p_j\|\geq1$ for every 
$j\in\mathbb N$. Consequently $\ell^1(w)\subseteq\ell^1$; in particular, every sequence in $\ell^1(w)$ converges 
to $0$. If $(\lambda_j)_{j\in\mathbb N}\in\ell^1(w)$, then 
$\sum\limits_{j=1}^{\infty}\|\lambda_jp_j\|=\sum\limits_{j=1}^{\infty}|\lambda_j|w_j<\infty$ and therefore
the series converges absolutely, which gives
$$a=\sum_{j=1}^{\infty}\lambda_jp_j\in \mathscr A.$$
By pairwise orthogonality, 
\begin{equation}\label{eq:a-acts-on-pj}
ap_j=p_ja=\lambda_jp_j\qquad (j\in\mathbb N).
\end{equation}

We first prove that $\{0\}\cup\{\lambda_j;\;j\in\mathbb N\}\subseteq\sigma(a)$. Fix $j\in\mathbb N$. By 
\eqref{eq:a-acts-on-pj} we have $(a-\lambda_je)p_j=0$. If $a-\lambda_je$ were invertible with inverse $c$, 
then $p_j=c(a-\lambda_je)p_j=0$, a contradiction. Hence, $\lambda_j\in\sigma(a)$ for every $j\in\mathbb N$. 
Since $\lambda_j\to0$ and $\sigma(a)$ is closed, also $0\in\sigma(a)$.
Conversely, let $\mu\notin\{0\}\cup\{\lambda_j;\;j\in\mathbb N\}$. The set 
$\{0\}\cup\{\lambda_j;\;j\in\mathbb N\}$ is compact since it is the union of a convergent sequence and its limit.
Let $\delta=\dist\big(\mu,\{0\}\cup\{\lambda_j;\;j\in\mathbb N\}\big)>0$. Hence, 
$|\lambda_j-\mu|\geq\delta$ for all $j\in\mathbb N$ and $|\mu|\geq\delta$. Put
$\gamma_j=\lambda_j \mu^{-1}(\lambda_j-\mu)^{-1}$ $(j\in\mathbb N)$.
Then $|\gamma_j|\leq\delta^{-2}|\lambda_j|$, whence $\sum\limits_{j=1}^{\infty}|\gamma_j|w_j\leq 
\delta^{-2}\sum\limits_{j=1}^{\infty}|\lambda_j|w_j<\infty$. It follows that
$$b=\sum_{j=1}^\infty\gamma_jp_j-\mu^{-1}e$$
is a well-defined element of $\mathscr A$. Using \eqref{eq:a-acts-on-pj} and pairwise orthogonality, we obtain
\begin{align*}
(a-\mu e)b&=a\sum_{j=1}^{\infty}\gamma_j p_j-\mu^{-1}a-\mu\sum_{j=1}^{\infty}\gamma_j p_j+e\\
&=\sum_{j=1}^{\infty}\big(\gamma_j\lambda_j-\mu^{-1}\lambda_j-\mu\gamma_j\big)p_j+e=e,
\end{align*}
because $\gamma_j\lambda_j-\mu^{-1}\lambda_j-\mu\gamma_j=0$ for every $j\in\mathbb N$. 
Since $a$ and $b$ commute, we see that $b$ is the inverse of $a-\mu e$, which means that $\mu\not\in\sigma(a)$. 
This proves $\sigma(a)=\{0\}\cup\{\lambda_j;\;j\in\mathbb N\}$.

Assume finally that the numbers $\lambda_j$ are pairwise distinct and nonzero. Since $\lambda_j\to0$, each 
$\lambda_j$ is an isolated point of $\sigma(a)$. It follows that the singletons $\{\lambda_j\}$ $(j\in\mathbb N)$ are 
pairwise disjoint components of $\sigma(a)$. Hence, $\sigma(a)$ has infinitely many components, that is, 
$a\in\mathcal C_\infty(\mathscr A)$.
\end{proof}

\begin{corollary}
\label{cor:Cinfty-infinite-orthogonal-family}
If $\mathscr A$ contains a sequence of nonzero pairwise orthogonal idempotents, then 
$\mathcal C_\infty(\mathscr A)\neq\varnothing$.
\end{corollary}

\begin{proof}
Let $(p_j)_{j\in\mathbb N}$ be such a sequence of idempotents. Let $w_j=\|p_j\|$ and define recursively 
$\lambda_0=1$ and $\lambda_j=\min\{2^{-j}w_j^{-1},\,\lambda_{j-1}/2\}$ for $j\in\mathbb N$. 
It is easily seen that $(\lambda_j)_{j\in\mathbb N}\in\ell^1(w)$ is a sequence of pairwise distinct and nonzero
numbers. Now apply Theorem~\ref{thm:Cinfty-orthogonal-idempotents}.
\end{proof}

We now determine the spectral depth of finite-dimensional algebras. If $\mathscr A$ is a finite-dimensional unital
complex algebra, then it can be embedded into the algebra $\mathbb M_n$ of all complex $n\times n$ matrices
for some $n\in\mathbb N$ such that the unit of $\mathscr A$ is mapped into the identity matrix $I$
(see \cite[Proposition 2.38]{Bre14} and its proof). Hence, there is no loss of generality if we assume that
$\mathscr A$ is a subalgebra of $\mathbb M_n$ and $I$ is the unit of $\mathscr A$. Let $\mathscr N$ be the
radical of $\mathscr A$, that is, the maximal nilpotent ideal. Then the quotient $\mathscr A'=\mathscr A/\mathscr N$
is a semisimple algebra. By the Wedderburn--Artin theorem, there exist positive integers $n_1,\ldots,n_r$
such that 
\begin{equation} \label{eq:wedderburn--artin}
\mathscr A'\cong\mathbb M_{n_1}\oplus\cdots\oplus\mathbb M_{n_r}\quad\text{ and }\quad n_1+\cdots+n_r\leq n.
\end{equation}
The reader is referred to Sections 2.9 and 2.10 in \cite{Bre14} for details.

\begin{proposition}
\label{prop:depth-finite-dim-algebra}
Let $\mathscr A$ and $\mathscr A'$ be as above and let \eqref{eq:wedderburn--artin} hold. Then
$$d_\sigma(\mathscr A) = d_\sigma(\mathscr A') = n_1+\cdots+n_r.$$
\end{proposition}

\begin{proof}
Let $\pi\colon \mathscr A \to \mathscr A'=\mathscr A/\mathscr N$ be the quotient map. We claim that an 
element $A\in\mathscr A$ is invertible if and only if $\pi(A)$ is invertible in
$\mathscr A'$. Indeed, suppose that $A\in\mathscr A$ and $\pi(A)$ is invertible in $\mathscr A'$.
It follows that there exist $B\in\mathscr A$ and $N\in\mathscr N$ such that $AB=I+N$. 
Since $\mathscr N$ is a nilpotent ideal, $N$ is nilpotent and therefore $I+N$ is invertible. This implies that $A$ 
has a right inverse. Similarly, there exist $B'\in\mathscr A$ and $N'\in\mathscr N$ such that $B'A=I+N'$. Since
$I+N'$ is invertible, $A$ has a left inverse. We conclude that $A$ is invertible in 
$\mathscr A$. The opposite is obvious: if $A$ is invertible, then $\pi(A)$ is invertible.
Applying the claim to $\lambda I-A$, we obtain $\sigma_{\mathscr A}(A)=\sigma_{\mathscr A'}(\pi(A))$ for
every $A\in\mathscr A$.
Since $\pi$ is a surjective mapping, it follows that $\mathscr A$ and $\mathscr A'$ have exactly the same set of 
achievable spectra, and therefore $d_\sigma(\mathscr A) = d_\sigma(\mathscr A')$.

Under the identification $\mathscr A' \cong \mathbb{M}_{n_1}\oplus \cdots \oplus \mathbb{M}_{n_r}$, an element 
$T=T_1\oplus\cdots\oplus T_r \in\mathscr A'$ is invertible if and only if each $T_i$ is invertible in 
$\mathbb{M}_{n_i}$. Consequently,
$$ \sigma_{\mathscr A'}(T) = \bigcup_{i=1}^r \sigma(T_i).$$
Each spectrum $\sigma(T_i)$ is a finite set of eigenvalues of $T_i$, of cardinality at most $n_i$. It follows that
$\sigma_{\mathscr A'}(T)$ is finite with cardinality at most $n_1+\cdots+n_r$.
As $\sigma_{\mathscr A'}(T)$ is a finite subset of $\mathbb{C}$, its connected components coincide with its points. 
Hence, the number of spectral components of $T$ equals the cardinality of $\sigma_{\mathscr A'}(T)$.

To see that the bound $n_1+\cdots+n_r$ is attained, choose each $T_i \in \mathbb{M}_{n_i}$ to be a diagonal 
matrix with $n_i$ distinct entries, and choose these entries, across all $i=1,\ldots,r$, to be pairwise distinct 
complex numbers. For this specific $T$, $\sigma_{\mathscr A'}(T)$ consists of exactly $n_1+\cdots+n_r$ 
distinct points. Therefore, $d_\sigma(\mathscr A') = n_1+\cdots+n_r$; the proposition follows.
\end{proof}

\begin{corollary} \label{cor:depth-finite-dim-Banach-space}
If $X$ is a complex Banach space of finite dimension $n\geq 2$, then $d_\sigma(\mathscr B(X))=n$.
\end{corollary}

\begin{proof}
The algebra $\mathscr B(X)$ can be identified with the matrix algebra $\mathbb{M}_n$, which is a simple
algebra and therefore $r=1$ and $n_1=n$. Hence, by Proposition \ref{prop:depth-finite-dim-algebra}, the 
spectral depth is $n$.
\end{proof}

It is not hard to see that for every integer $1\leq r\leq n$, there exists a subalgebra 
$\mathscr A\subseteq \mathscr B(X)$ with $I\in \mathscr A$ such that $d_\sigma(\mathscr A)=r$.
Indeed, choose integers $n_1,\ldots, n_r$ such that $n=n_1+\cdots+n_r$ and let $X=X_1\oplus\cdots\oplus X_r$ 
be a decomposition with $\dim(X_j)=n_j$ for $j=1,\ldots, r$.
Denote by $I_j$ the identity operator on $X_j$ and let
$$\mathscr A=\{ \alpha_1 I_1\oplus\cdots\oplus \alpha_r I_r\in\mathscr B(X);\; \alpha_1,\ldots,
\alpha_r\in\mathbb{C}\}.$$
It is easily seen that $\mathscr A$ is a subalgebra of $\mathscr B(X)$ containing $I$. The spectrum of
$A=\alpha_1 I_1\oplus\cdots\oplus \alpha_r I_r$ is $\sigma(A)=\{\alpha_1,\ldots,\alpha_r\}$. Hence,
the cardinality of the spectrum of an operator in $\mathscr A$ is at most $r$; it is precisely $r$ if $\alpha_j$
are distinct numbers. Thus, $d_\sigma(\mathscr A)=r$.\medskip

We now turn to infinite-dimensional Banach spaces. The following is a corollary of Theorem 
\ref{thm:Cinfty-orthogonal-idempotents}.

\begin{corollary}
\label{cor:spectral-depth-BX}
If $X$ is an infinite-dimensional complex Banach space, then for every sequence 
$(\lambda_j)_{j\in\mathbb N}\in\ell^1$ there exists $T\in \mathscr K(X)$ such that 
$\sigma(T)=\{0\}\cup\{\lambda_j;\; j\in\mathbb N\}$. In particular, if the 
$\lambda_j$ are pairwise distinct and nonzero, then $T\in \mathcal C_\infty(\mathscr B(X))$ and, consequently, 
$d_\sigma(\mathscr B(X))=\infty$.
\end{corollary}

\begin{proof}
Every infinite-dimensional Banach space contains a normalized basic sequence 
$(e_j)_{j\in\mathbb N}\subseteq X$ with basis constant $K$ and the corresponding biorthogonal functionals 
$(\psi_j)_{j\in\mathbb N}\subseteq Y^*$, where $Y=\overline{\Span}\{e_j;\;j\in\mathbb N\}$, satisfy $\|\psi_j\|\leq 2K$; 
see \cite[Section~1.5]{AK06}. For each $j\in\mathbb N$, let $\xi_j\in X^*$ be a Hahn--Banach extension of $\psi_j$, 
so that $\langle e_i,\xi_j\rangle=\delta_{ij}$ and $\|\xi_j\|\leq 2K$. Put $P_j=e_j\otimes \xi_j$. Then 
$$P_j^2=\langle e_j,\xi_j\rangle\, e_j\otimes \xi_j=P_j \qquad \textrm{and}\qquad  
P_iP_j=\langle e_j,\xi_i\rangle\, e_i\otimes \xi_j=0\quad (i\neq j),$$
so $(P_j)_{j\in\mathbb N}$ is a sequence of nonzero pairwise orthogonal idempotents in $\mathscr B(X)$ with 
$\|P_j\|=\|e_j\|\,\|\xi_j\|\leq 2K$. Hence, $\ell^1\subseteq\ell^1(w)$ for $w_j=\|P_j\|$, and 
Theorem~\ref{thm:Cinfty-orthogonal-idempotents} yields, for every $(\lambda_j)_{j\in\mathbb N}\in\ell^1$, an 
operator $T=\sum_{j}\lambda_jP_j$ with $\sigma(T)=\{0\}\cup\{\lambda_j;\;j\in\mathbb N\}$. The operator $T$ is 
compact, being the norm limit of the finite-rank partial sums $\sum_{j\leq n}\lambda_jP_j$. The last assertion 
is the final statement of Theorem~\ref{thm:Cinfty-orthogonal-idempotents}.
\end{proof}

Theorem~\ref{thm:Cinfty-orthogonal-idempotents} applies verbatim to subalgebras of $\mathscr B(X)$, and it does 
so for an arbitrary complete algebra norm, not only for the operator norm. Recall that a subalgebra of 
$\mathscr B(X)$ containing $\mathscr F(X)$ is often called a \emph{standard operator algebra}.

\begin{corollary}
\label{cor:Cinfty-standard-operator-algebras}
Let $X$ be an infinite-dimensional complex Banach space and let $\mathscr A$ be a subalgebra of $\mathscr B(X)$ 
which contains $I$ and $\mathscr F(X)$, and which is a Banach algebra with respect to some norm $\|\cdot\|'$. 
Then there exists an operator $T$ in the $\|\cdot\|'$-closure of $\mathscr F(X)$ in $\mathscr A$ such that 
$T\in\mathcal C_\infty(\mathscr A)$; in particular, $d_\sigma(\mathscr A)=\infty$. If, moreover, the idempotents 
$P_j=e_j\otimes \xi_j$ constructed in the proof of Corollary~\ref{cor:spectral-depth-BX} satisfy 
$\sup_{j\in\mathbb N}\|P_j\|'<\infty$, then the spectrum of $T$ in $\mathscr A$ can be prescribed to be 
$\{0\}\cup\{\lambda_j;\;j\in\mathbb N\}$ for an arbitrary sequence $(\lambda_j)_{j\in\mathbb N}\in\ell^1$ of 
pairwise distinct nonzero numbers.
\end{corollary}

\begin{proof}
The idempotents $P_j$ of the proof of Corollary~\ref{cor:spectral-depth-BX} are of rank one, hence they 
belong to $\mathscr F(X)\subseteq\mathscr A$, and they are nonzero and pairwise orthogonal. Applying 
Corollary~\ref{cor:Cinfty-infinite-orthogonal-family} to the Banach algebra $(\mathscr A,\|\cdot\|')$ we obtain 
$T=\sum_j\lambda_jP_j\in\mathcal C_\infty(\mathscr A)$, and $T$ is the $\|\cdot\|'$-limit of the finite-rank 
operators $\sum_{j\leq n}\lambda_jP_j$. The second assertion follows from 
Theorem~\ref{thm:Cinfty-orthogonal-idempotents}, because $\sup_j\|P_j\|'<\infty$ gives 
$\ell^1\subseteq\ell^1(w)$ for $w_j=\|P_j\|'$.
\end{proof}

In Corollary~\ref{cor:Cinfty-standard-operator-algebras} the spectrum is computed in $\mathscr A$ and not in 
$\mathscr B(X)$ because the two may differ, since $\mathscr A$ need not be inverse-closed in $\mathscr B(X)$. 
The proof of Theorem~\ref{thm:Cinfty-orthogonal-idempotents} avoids this difficulty, because the inverse of
$a-\mu e$ is exhibited explicitly as an element of the algebra under consideration. 

The boundedness hypothesis in the second assertion of Corollary \ref{cor:Cinfty-standard-operator-algebras}
is satisfied in the natural examples. For instance, let $H$ be an infinite-dimensional complex Hilbert space and
let $\mathscr S_p(H)\subseteq \mathscr B(H)$ be the Schatten $p$-class $(1\leq p<\infty)$. 
Note that $\mathscr S_p(H)$ contains all finite-rank operators. Since $\mathscr S_p(H)$ is a proper ideal, it does 
not contain $I$ and therefore $\mathscr A=\mathbb CI+\mathscr S_p(H)$ is a well-defined subalgebra of 
$\mathscr B(H)$ containing $I$. Since $\mathscr S_p(H)$ is closed in the norm $\|\cdot\|_p\geq \|\cdot\|$, one
can define a norm $\|\alpha I+S\|'=|\alpha|+\|S\|_p$ on $\mathscr A$. One can show that this norm
is submultiplicative. Thus, $\mathscr A$ is a Banach algebra containing $I$ and all finite-rank operators. The norm 
of every rank-one operator is $\|x\otimes \xi\|_p=\|x\|\,\|\xi\|$; in particular $\sup_j\|P_j\|'<\infty$. Hence,
this algebra satisfies the assumptions of Corollary \ref{cor:Cinfty-standard-operator-algebras}.

\section{Topological properties of spectral classes} \label{sec:topological-properties-spectral-classes}
In this section we investigate the topological structure of the spectral classes $\mathcal C_n(\mathscr A)$. 
We will need the following continuity property of the Riesz functional calculus in the element, for a fixed 
holomorphic function on a fixed open set. We include a proof for the reader's convenience.

\begin{lemma}
\label{lem:continuity-functional-calculus}
Let $a\in\mathscr A$ and let $f$ be a holomorphic function on a neighbourhood $U$ of $\sigma(a)$. For every 
$\varepsilon>0$ there exists $\delta>0$ such that $\sigma(b)\subseteq U$ and $\|f(a)-f(b)\|<\varepsilon$
for every $b\in\mathscr A$ satisfying $\| a-b\|<\delta$.
\end{lemma}

\begin{proof}
If $f$ is identically equal to $0$, the statement is trivial. Suppose therefore that $f\ne 0$.
Fix a finite collection $\Gamma$ of positively oriented Jordan curves contained in $U\setminus\sigma(a)$ and 
enclosing $\sigma(a)$. Let $l(\Gamma)$ be the total length of $\Gamma$. 
Hence, $\sigma(a)\subseteq\ins(\Gamma)$ and $\Gamma\cap\ins(\Gamma)=\varnothing$. It follows that 
$\Gamma$ is a compact subset of the resolvent set $\rho(a)$ and therefore the norm of the resolvent is a continuous 
function on it, achieving the maximal value 
$$M=\max\{\|(\lambda e-a)^{-1}\|;\;\lambda\in\Gamma\}.$$
Denote $M_f=\max\{|f(\lambda)|;\;\lambda\in\Gamma\}$; since $f\ne 0$, this number is positive.

By the upper semicontinuity of the spectrum \cite[Theorem~1]{New51} there exists $\delta_0>0$ such that 
$\sigma(b)\subseteq\ins(\Gamma)$ whenever $\|b-a\|<\delta_0$. This means that $\sigma(b)\subseteq U$ and 
the same collection $\Gamma$ encloses $\sigma(b)$. Hence, $f(b)$ is defined by \eqref{eq:holomorphic-calculus}.
Let
$$\delta=\min\Big\{\delta_0,\ \frac{1}{2M},\ \frac{\pi\varepsilon}{M^2 l(\Gamma)M_f}\Big\}.$$

If $b\in\mathscr A$ is such that $\|b-a\|<\delta$ and $\lambda\in\Gamma$, then
$$\lambda e-b=\lambda e-a-(b-a)=(\lambda e-a)\big(e-(\lambda e-a)^{-1}(b-a)\big).$$
Since $\|(\lambda e-a)^{-1}(b-a)\|\leq\|(\lambda e-a)^{-1}\|\|b-a\|\leq M\|b-a\|<M\delta\leq\frac{1}{2}$, 
we see that $e-(\lambda e-a)^{-1}(b-a)$ is invertible.
It follows that $\lambda e-b$ is invertible. Moreover, since
$$\big\|\big(e-(\lambda e-a)^{-1}(b-a)\big)^{-1}\big\|\leq
\sum_{j=0}^{\infty}\big\|(\lambda e-a)^{-1}(b-a)\big\|^j
<\sum_{j=0}^{\infty} 2^{-j}=2,$$
we have $\|(\lambda e-b)^{-1}\|< 2M$. The resolvent identity
$$(\lambda e-b)^{-1}-(\lambda e-a)^{-1}=(\lambda e-b)^{-1}(b-a)(\lambda e-a)^{-1}$$
then gives $\|(\lambda e-b)^{-1}-(\lambda e-a)^{-1}\|\leq 2M^2\|b-a\|$ for every $\lambda\in\Gamma$.

The elements $f(a)$ and $f(b)$ are both given by integrals over $\Gamma$, therefore
$$\|f(a)-f(b)\|=\Big\|\frac{1}{2\pi i}\int_\Gamma f(\lambda)
\big((\lambda e-a)^{-1}-(\lambda e-b)^{-1}\big)\,d\lambda\Big\|
\leq\frac{M^2 l(\Gamma)M_f}{\pi}\,\|b-a\|<\varepsilon$$
by the choice of $\delta$.
\end{proof}

\begin{theorem}
\label{thm:Ck-open}
For every $n\in\mathbb N$, the spectral class $\mathcal C_n(\mathscr A)$ is an open subset of $\mathscr A$.
\end{theorem}

\begin{proof}
The case $n=1$ is clear, since $\mathcal C_1(\mathscr A)=\mathscr A$. Let therefore $n\geq2$ and let 
$a\in\mathcal C_n(\mathscr A)$. As $\sigma(a)$ has at least $n$ components, 
Lemma~\ref{lem:finite-clopen-spectral-partition} provides nonempty pairwise disjoint clopen subsets 
$\sigma_1,\ldots,\sigma_n$ of $\sigma(a)$ with $\sigma(a)=\sigma_1\cup\cdots\cup\sigma_n$. Being compact and 
pairwise disjoint, these sets admit pairwise disjoint open neighbourhoods $U_1,\ldots,U_n$; put 
$U=U_1\cup\cdots\cup U_n$.

For $j=1,\ldots,n$ let $\chi_j\colon U\to\mathbb C$ be the characteristic function of $U_j$. Since 
$U_1,\ldots,U_n$ are open and pairwise disjoint, $\chi_j$ is locally constant, hence holomorphic on the 
neighbourhood $U$ of $\sigma(a)$. Since $\sigma(a)\cap U_j=\sigma_j$, we see that $\chi_j(a)=e_{\sigma_j}$,
the spectral idempotent of $a$ associated with the nonempty clopen subset $\sigma_j$; hence $\chi_j(a)\ne 0$ by 
Proposition~\ref{prop:orthogonal-idempotent-decomposition}. Applying 
Lemma~\ref{lem:continuity-functional-calculus} to $a$, $U$, $\chi_j$ and $\varepsilon_j=\|\chi_j(a)\|$, we obtain 
$\delta_j>0$ such that
$$\sigma(b)\subseteq U\quad\text{and}\quad \|\chi_j(a)-\chi_j(b)\|<\|\chi_j(a)\|
\qquad\text{whenever } \|a-b\|<\delta_j .$$

Put $\delta=\min\{\delta_1,\ldots,\delta_n\}$ and let $b\in\mathscr A$ satisfy $\|a-b\|<\delta$. Fix 
$j\in\{1,\ldots,n\}$. Then $\chi_j(b)\neq0$, since otherwise the displayed inequality would read 
$\|\chi_j(a)\|<\|\chi_j(a)\|$. On the other hand, since $\sigma(b)\subseteq U$, the set $\sigma(b)\cap U_j$ is 
a clopen subset of $\sigma(b)$ and $\chi_j(b)$ is precisely the associated spectral idempotent. Were 
$\sigma(b)\cap U_j$ empty, the function $\chi_j$ would vanish on a neighbourhood of $\sigma(b)$ and we would 
get $\chi_j(b)=0$ contradicting $\chi_j(b)\ne 0$. Hence, $\sigma(b)\cap U_j\neq\varnothing$ for every $j=1,\ldots,n$.

The sets $\sigma(b)\cap U_1,\ldots,\sigma(b)\cap U_n$ are thus nonempty, pairwise disjoint, and relatively 
clopen in $\sigma(b)$, and their union is $\sigma(b)$. Since a connected subset of $\sigma(b)$ meeting a 
relatively clopen set is contained in it, each component of $\sigma(b)$ lies in exactly one of them, and each 
of them contains at least one component. Therefore, $\sigma(b)$ has at least $n$ components, that is, 
$b\in\mathcal C_n(\mathscr A)$. This proves that $\mathcal C_n(\mathscr A)$ is open.
\end{proof}

Let $\mathscr I$ be a proper closed two-sided ideal in $\mathscr A$ and let 
$\pi\colon\mathscr A\to\mathscr A/\mathscr I$ be the quotient map. Then $\mathscr A/\mathscr I$ is a complex 
Banach algebra with unit $\pi(e)$. Recall that $\pi$ is a continuous, open and surjective mapping. We shall use 
the following elementary fact.

\begin{lemma}\label{lem:open-quotient-maps-density-both-ways}
Let $\phi\colon\Omega_1\to\Omega_2$ be a continuous, open and surjective mapping between topological spaces. 
Then a subset $\Delta\subseteq\Omega_2$ is dense in $\Omega_2$ if and only if the set $\phi^{-1}(\Delta)$ is 
dense in $\Omega_1$.
\end{lemma}

\begin{proof}
Assume that $\Delta$ is dense in $\Omega_2$ and let $U\subseteq\Omega_1$ be a nonempty open set. 
Then $\phi(U)$ is nonempty and open and therefore $\phi(U)\cap \Delta\neq\varnothing$.
This gives $U\cap \phi^{-1}(\Delta)\neq\varnothing$. Hence, $\phi^{-1}(\Delta)$ is dense in $\Omega_1$. 
Conversely, assume that $\phi^{-1}(\Delta)$ is dense in $\Omega_1$. Let $V\subseteq\Omega_2$ be a nonempty 
open set. By surjectivity, $\phi^{-1}(V)$ is nonempty, and it is open by continuity. Hence, 
$\phi^{-1}(V)\cap \phi^{-1}(\Delta)\neq\varnothing$, and therefore $V\cap \Delta\neq\varnothing$.
\end{proof}

\begin{proposition}
\label{prop:quotient-Ck-basic-properties}
Let $\mathscr I$ be a proper closed two-sided ideal in $\mathscr A$ and let $n\in\mathbb N$. Then the set 
$\pi^{-1}\big(\mathcal C_n(\mathscr A/\mathscr I)\big)$ is open in $\mathscr A$. 
If $a\in\pi^{-1}\big(\mathcal C_n(\mathscr A/\mathscr I\big)$, then 
$a+b\in \pi^{-1}\big(\mathcal C_n(\mathscr A/\mathscr I)\big)$ for all $b\in\mathscr I$.
Moreover, $\pi^{-1}\big(\mathcal C_n(\mathscr A/\mathscr I)\big)$ is dense in $\mathscr A$ if and only if 
$\mathcal C_n(\mathscr A/\mathscr I)$ is dense in $\mathscr A/\mathscr I$.
\end{proposition}

\begin{proof}
The class $\mathcal C_n(\mathscr A/\mathscr I)$ is open in $\mathscr A/\mathscr I$ by 
Theorem~\ref{thm:Ck-open}, and $\pi$ is continuous; hence its preimage is open in $\mathscr A$.
If $b\in\mathscr I$, then $\pi(a+b)=\pi(a)$, which gives the second assertion. The last assertion follows 
from Lemma~\ref{lem:open-quotient-maps-density-both-ways}.
\end{proof}

The openness established in Theorem \ref{thm:Ck-open} holds for every complex Banach algebra. In contrast, 
as the following results show, the existence and abundance of elements in the spectral classes 
$\mathcal C_n(\mathscr A)$ depends strongly on the algebra $\mathscr A$. 

\begin{proposition}
\label{prop:Ck-spectral-depth}
Let $\Omega$ be a compact Hausdorff space and let $C(\Omega)$ be the Banach algebra of all continuous complex 
functions defined on $\Omega$. The following statements hold. 
\begin{enumerate}
\item $\mathcal C_n\big(C(\Omega)\big)\ne\varnothing$ if and only if $\Omega$ has at least $n$ components.
\item $\mathcal C_\infty\big(C(\Omega)\big)\ne\varnothing$ if and only if $\Omega$ has infinitely many components.
\end{enumerate}
\end{proposition}

\begin{proof}
(i) Suppose first that $\Omega=\Delta_1\cup\cdots\cup \Delta_m$ is the decomposition of $\Omega$ into its $m$ 
components, where $m<n$. For $f\in C(\Omega)$, each set $f(\Delta_j)$ is connected, and therefore 
$\sigma(f)=f(\Omega)=\bigcup_{j=1}^m f(\Delta_j)$ has at most $m$ components. This proves that
$\mathcal C_n\big(C(\Omega)\big)=\varnothing$.
Conversely, if $\Omega$ has at least $n$ components, apply Lemma~\ref{lem:finite-clopen-spectral-partition} 
to obtain a clopen partition $\Omega=\sigma_1 \cup \cdots\cup \sigma_n$ with $\sigma_j$ nonempty for each 
$1\leq j\leq n$. For distinct scalars $\lambda_1,\ldots,\lambda_n$, the function
$f=\sum_{j=1}^n\lambda_j\chi_{\sigma_j}$ has spectrum $\sigma(f)=f(\Omega)=\{\lambda_1,\ldots,\lambda_n\}$, 
which means that $f\in \mathcal C_n\big(C(\Omega)\big)$.

(ii) Suppose that $\Omega$ has infinitely many components $(\Delta_j)_{j\in\mathbb N}$. 
By \cite[Theorem~6.1.23]{Eng89} there exists a nonempty clopen set $\sigma$ such that $\Delta_1\subseteq \sigma$ 
and $\Delta_2\cap \sigma=\varnothing$. Since each connected set is either contained in a given clopen set or its 
complement, one of the clopen sets $\sigma$ and $\Omega\setminus \sigma$ contains infinitely many components 
among $(\Delta_j)_{j\in\mathbb N}$. Let us denote this clopen set by $\sigma^{(1)}$. Taking 
$\sigma_1=\Omega\setminus \sigma^{(1)}$ (nonempty and clopen) and repeating the argument within 
$\sigma^{(1)}$ we obtain $\sigma^{(2)}\subseteq \sigma^{(1)}$ with infinitely many components, and we set 
$\sigma_2=\sigma^{(1)}\setminus\sigma^{(2)}$. Continuing in this way, we obtain a sequence 
$(\sigma_j)_{j\in\mathbb N}$ of nonempty pairwise disjoint clopen sets. The series 
$\sum_{j=1}^{\infty}2^{-j}\chi_{\sigma_j}$ converges uniformly; let $f\in C(\Omega)$ denote its sum. Set 
$F=\Omega\setminus \bigcup_{j=1}^\infty \sigma_j$. For every $k\in\mathbb N$, the set 
$\Omega\setminus \big(\sigma_1\cup\cdots\cup\sigma_k\big)=\sigma^{(k)}$ is nonempty and compact. It follows
from $F=\bigcap_{k=1}^{\infty}\sigma^{(k)}$ that $F$ is nonempty, by the finite intersection property of 
compact spaces. Since $f(x)=2^{-j}$ for each $x\in \sigma_j$ and $f(x)=0$ for each $x\in F$, we conclude 
$$\sigma(f)=f(\Omega)=\{0\}\cup\{2^{-j};\; j\in\mathbb N\},$$
that is, $\sigma(f)$ has infinitely many components. Conversely, if $\Omega$ has only finitely many components,
say $m$, then by the argument in (i) every $f\in C(\Omega)$ has spectrum with at most $m$ components.
This means that $\mathcal C_{m+1}\big(C(\Omega)\big)=\varnothing$ and hence
$\mathcal C_{\infty}\big(C(\Omega)\big)=\varnothing$.
\end{proof}

Recall that a topological space $\Omega$ is totally disconnected if the components of $\Omega$ are one-point
sets (e.g., the Cantor set is totally disconnected).
 
\begin{proposition}
\label{prop:Cn-dense-totally-disconnected}
Let $\Omega$ be a totally disconnected compact Hausdorff space with infinitely many points. Then 
$\mathcal C_n\big(C(\Omega)\big)$ is dense in $C(\Omega)$ for every $n\in\mathbb N$.
\end{proposition}

\begin{proof}
Since $\Omega$ is compact Hausdorff and totally disconnected, the component of each point is a singleton, and 
by \cite[Theorem~6.1.23]{Eng89} it is the intersection of all clopen sets containing that point. A standard 
compactness argument then shows that $\Omega$ has a base of clopen sets. 
This implies that locally constant functions, i.e., functions which are of the form $\sum_{j=1}^{k}\mu_j\chi_{\Delta_j}$, where $\mu_1,\ldots,\mu_k\in\mathbb C$ and $\Delta_1,\ldots,\Delta_k$ is a clopen partition of $\Omega$, 
form a self-adjoint algebra that contains the constant functions and separates points. 
By the Stone--Weierstrass theorem \cite[Theorem~V.8.1]{Con90}, this algebra is dense in $C(\Omega)$.

Let $n\geq 2$ (the case $n=1$ is trivial). Let $f\in C(\Omega)$ and $\varepsilon>0$ be arbitrary.
As observed in the previous paragraph, we can choose a locally constant function $g\in C(\Omega)$ 
with $\|f-g\|<\varepsilon/2$. Hence, there are a clopen partition $\Omega=\Delta_1\cup\cdots\cup\Delta_k$ into 
nonempty sets and numbers $\mu_1,\ldots,\mu_k$ such that $g=\sum_{j=1}^{k}\mu_j\chi_{\Delta_j}$. 
We may assume that $k\geq n$. Indeed, since $\Omega$ is infinite, some $\Delta_j$ contains two distinct points 
$x$ and $y$. Because the space $\Omega$ has a base of clopen sets, there is a clopen set 
$\Delta_{j}^{\prime}\subseteq\Omega$ with $x\in \Delta_{j}^{\prime}$ and $y\notin \Delta_{j}^{\prime}$. 
Replacing $\Delta_j$ by the two nonempty clopen sets $\Delta_j\cap \Delta_{j}^{\prime}$ and 
$\Delta_j\setminus \Delta_{j}^{\prime}$, on each of which $g$ is still constant, we obtain a finer partition with 
one more member. Repeating this finitely many times, we may assume $k\geq n$.

Choose pairwise distinct numbers $\nu_1,\ldots,\nu_k$ with $|\nu_j-\mu_j|<\varepsilon/2$ for each $j=1,\ldots,k$. 
Let $h=\sum_{j=1}^{k}\nu_j\chi_{\Delta_j}\in C(\Omega)$. Then $\|g-h\|<\varepsilon/2$, and 
therefore $\|f-h\|<\varepsilon$. Since the sets $\Delta_j$ are nonempty and the numbers $\nu_j$ are pairwise 
distinct, we have $\sigma(h)=h(\Omega)=\{\nu_1,\ldots,\nu_k\}$, a set with $k\geq n$ components. Hence, 
$h\in\mathcal C_n\big(C(\Omega)\big)$, and the proof is complete.
\end{proof}

Note that the condition that $\Omega$ has infinitely many points is necessary in Proposition 
\ref{prop:Cn-dense-totally-disconnected}. Indeed, if $\Omega$ is finite with cardinality $k\in\mathbb N$,
then $\mathcal C_n\big(C(\Omega)\big)=\varnothing$ for every $n>k$, by 
Proposition~\ref{prop:Ck-spectral-depth}~(i). Now we prove that $\mathcal C_n(\mathscr B(X))$ is dense in 
$\mathscr B(X)$ for each $n\in \mathbb N$, if $X$ is infinite-dimensional.

\begin{theorem}\label{thm:Cn-density}
Let $X$ be an infinite-dimensional complex Banach space. Then for every $n\geq 2$, for every $\varepsilon>0$ and 
for every $T\in\mathscr B(X)$ there exists an operator $F\in\mathscr F_{n-1}(X)$ such that
$\|F\|<\varepsilon$ and $T+F\in\mathcal C_n\big(\mathscr B(X)\big)$.
\end{theorem}

\begin{proof}
The case $n=2$ is precisely the rank-one perturbation theorem of Herrero and Salinas \cite{HS72}. 
Assume, therefore, that $n\ge3$ and let $T\in\mathscr B(X)$ and $\varepsilon>0$ be arbitrary. 

If $\sigma(T)$ has at least $n$ components, then we can take $F=0$. Hence, we may assume that $\sigma(T)$ has 
$k<n$ components denoted by $\sigma_1,\ldots,\sigma_k$. For each $1\leq j\leq k$, let $P_j$ be the Riesz projection 
associated with the spectral set $\sigma_j$, and let $X_j=\ran(P_j)$. Since the projections $P_1,\ldots,P_k$ 
are pairwise orthogonal and their sum is the identity operator on $X$, we have $X=X_1\oplus\cdots\oplus X_k$,
where each $X_j$ is a closed subspace of $X$ invariant under $T$. The operator $T$ decomposes as 
$T=T_1\oplus\cdots\oplus T_k$ with $T_j\in\mathscr B(X_j)$ and $\sigma(T_j)=\sigma_j$ for $j=1,\ldots,k$. 

Since $X$ is infinite-dimensional, at least one of the spaces $X_1,\ldots,X_k$ is infinite-dimensional. Without loss 
of generality, we may assume that $X_k$ is infinite-dimensional. Choose an open neighbourhood $U$ of $\sigma_k$ 
such that $U\cap\sigma_j=\varnothing$ for $1\leq j\leq k-1$. Once again by Newburgh's theorem 
\cite[Theorem~1]{New51}, there exists $\eta>0$ such that $\sigma(S)\subseteq U$ whenever $S\in\mathscr B(X_k)$ 
satisfies $\|S-T_k\|<\eta$. Set $\delta=\min\left\{\frac{\varepsilon}{n},\eta\right\}$.
Applying the theorem of Herrero and Salinas \cite{HS72} to the operator $T_k$, we obtain a vector $u_k\in X_k$ 
and a functional $\xi_k\in X_k^*$ such that $\|u_k\otimes \xi_k\|<\delta$ and 
$T_k+u_k\otimes \xi_k\in\mathcal C_2\big(\mathscr B(X_k)\big)$. 
Since $T_k+u_k\otimes \xi_k\in\mathcal C_2\big(\mathscr B(X_k)\big)$, the spectrum of the operator 
$T_k+u_k\otimes \xi_k$ has at least two components $\sigma_k^{(1)}$ and $\sigma_k^{(2)}$. It follows from 
$\|u_k\otimes \xi_k\|<\delta\le\eta$ that $\sigma(T_k+u_k\otimes \xi_k)\subseteq U$ which means that
$\sigma_1,\ldots,\sigma_{k-1},\sigma_k^{(1)},\sigma_k^{(2)}$ are pairwise disjoint connected sets.

Define $u\in X$ and $\xi\in X^*$ by $u=0\oplus\cdots\oplus0\oplus u_k$ and $\xi(x_1\oplus\cdots\oplus x_k)=\xi_k(x_k)$.
Then $u\otimes \xi= 0\oplus\cdots\oplus0\oplus(u_k\otimes \xi_k)$ and therefore 
$\|u\otimes \xi\| = \|u_k\otimes \xi_k\| < \delta \le \frac{\varepsilon}{n}$. From 
$T+u\otimes \xi = T_1\oplus\cdots\oplus T_{k-1}\oplus (T_k+u_k\otimes \xi_k)$ and 
$\sigma(T+u\otimes \xi)=\sigma(T_k+u_k\otimes \xi_k)\cup \bigcup_{j=1}^{k-1} \sigma(T_j)=
\sigma(T_k+u_k\otimes \xi_k) \cup \bigcup_{j=1}^{k-1} \sigma_j$ it follows that $\sigma(T+u\otimes \xi)$ has at 
least $k+1$ components $\sigma_1,\ldots,\sigma_{k-1},\sigma_k^{(1)}$ and $\sigma_k^{(2)}$, proving 
$T+u\otimes \xi\in\mathcal C_{k+1}\big(\mathscr B(X)\big)$. 

If $k+1\ge n$, then taking $F=u\otimes \xi$ completes the proof. We may therefore assume that $k+1<n$. 
Since $T+u\otimes \xi\in\mathcal C_{k+1}\big(\mathscr B(X)\big)$, the spectrum of $T+u\otimes \xi$ admits 
a decomposition $\sigma(T+u\otimes \xi)= \sigma_{1}^{\prime}\cup \cdots\cup\sigma_{k+1}^{\prime}$,
where $\sigma_{1}^{\prime},\ldots,\sigma_{k+1}^{\prime}$ are nonempty pairwise disjoint clopen subsets of 
the spectrum $\sigma(T+u\otimes \xi)$. As before, $\dim(X)=\infty$ implies that at least one of the corresponding 
spectral subspaces is infinite-dimensional. We may therefore repeat the above construction on that spectral subspace. 
At each step, the number of components of the spectrum increases by at least one, while the rank-one perturbation 
added at that step has norm less than $\min\{\varepsilon/n,\eta'\}$, where $\eta'$ is the positive number
supplied by Newburgh's theorem for that step. Since the initial spectrum has $k$ components, after at 
most $n-k$ steps we obtain an operator whose spectrum has at least $n$ components. The total perturbation has 
rank at most $n-k\leq n-1$ and norm less than $(n-k)\varepsilon/n <\varepsilon$.
\end{proof}

In the infinite-dimensional case, the openness and density of the spectral classes 
$\mathcal C_n\big(\mathscr B(X)\big)$ allow us to apply the Baire category theorem.

\begin{theorem}\label{thm:Cinfty-residual}
Let $X$ be an infinite-dimensional complex Banach space. Then the set $\mathcal C_\infty\big(\mathscr B(X)\big)$ 
is a dense $G_\delta$-subset of $\mathscr B(X)$. 
\end{theorem}

\begin{proof}
By Theorems \ref{thm:Ck-open} and \ref{thm:Cn-density}, the set $\mathcal C_n\big(\mathscr B(X)\big)$ is open 
and dense in $\mathscr B(X)$ for every $n\in\mathbb N$. It follows, by the Baire category theorem, that 
$$\mathcal C_\infty\big(\mathscr B(X)\big)=\bigcap_{n=1}^\infty \mathcal C_n\big(\mathscr B(X)\big)$$
is a dense $G_\delta$-subset of $\mathscr B(X)$.
\end{proof}

By Theorem \ref{thm:Ck-open}, every spectral class $\mathcal C_n(\mathscr A)$ is an open subset of $\mathscr A$.
On the other hand, the spectral class $\mathcal C_\infty(\mathscr A)$ need not be open. Indeed, let $X$ be an 
infinite-dimensional complex Banach space and let $T=\sum_{j=1}^\infty\lambda_jP_j$ be the compact operator
constructed in the proof of Corollary~\ref{cor:spectral-depth-BX} with pairwise distinct and nonzero numbers 
$\lambda_j$ $(j\in\mathbb N)$. For $m\in\mathbb N$ put $T_m=\sum_{j=1}^m\lambda_jP_j$. Then $T_m$ is 
of finite rank and therefore $T_m\notin\mathcal C_\infty\big(\mathscr B(X)\big)$. On the other hand, 
$\lim\limits_{m\to\infty}\|T-T_m\|=0$. Hence, no open ball centred at $T$ is contained in 
$\mathcal C_\infty\big(\mathscr B(X)\big)$.

We conclude this section with an example of a complex Banach algebra $\mathscr A$ with a unit 
whose spectral classes $\mathcal C_n(\mathscr A)$ $(n\geq 2)$ are nonempty and not dense in $\mathscr A$. 
 
\begin{example}\label{ex:all-spectral-classes-nonempty-none-dense}
Let $m\geq 2$ and let $\mathbb I$ be either the finite set $\{1,\ldots,m\}$ or $\mathbb N$ equipped with the 
discrete topology. Set $\Omega=\mathbb I\times \overline{\mathbb D}$, where $\mathbb D$ is the open unit disk 
in the complex plane. We will consider the commutative Banach algebra $C_b(\Omega)$ of all bounded continuous 
complex functions on $\Omega$. Then each element of this algebra can be seen as a finite or infinite sequence
$f=(f_j)_{j\in\mathbb I}$ of continuous complex functions on $\overline{\mathbb D}$; the norm of
$f$ is $\|f\|=\sup_{j\in\mathbb I}\|f_j\|$, where $\|f_j\|$ means the supremum norm on $\overline{\mathbb D}$.

If $\mathbb I=\{1,\ldots,m\}$, then $\Omega$ is a compact subset of $\mathbb C^2$ with components 
$\{j\}\times \overline{\mathbb D}$ for $1\leq j\leq m$. By Proposition~\ref{prop:Ck-spectral-depth} it follows 
that the spectral depth is $d_\sigma\big(C_b(\Omega)\big)=m$. 

Let us now consider the case $\mathbb I=\mathbb N$. We are going to prove that 
$\mathcal C_n\big(C_b(\Omega)\big)\neq\varnothing$ for every $n\in\mathbb N$, whereas 
$\mathcal C_n\big(C_b(\Omega)\big)$ is not dense in $C_b(\Omega)$ for any $n\geq2$. An easy exercise shows 
that the spectrum of an element $f=(f_j)_{j\in\mathbb N}\in C_b(\Omega)$ is 
$$\sigma(f)=\overline{\bigcup_{j\in\mathbb N}f_j(\overline{\mathbb D})}.$$
To prove that all spectral classes are nonempty, choose distinct nonzero complex numbers 
$(\lambda_j)_{j\in\mathbb N}$ converging to zero. For each $j\in\mathbb N$, let 
$c_{\lambda_j}\in C(\overline{\mathbb D})$ be the constant function equal to $\lambda_j$.
Define $f=(c_{\lambda_j})_{j\in\mathbb N}\in C_b(\Omega)$.
Then the description of the spectrum above yields $\sigma(f)=\{0\}\cup\{\lambda_j;\; j\in\mathbb N\}$ which 
shows that $f\in\mathcal C_n\big(C_b(\Omega)\big)$ for every $n\in\mathbb N$.

It remains to prove that for $n\geq 2$ the class $\mathcal C_n\big(C_b(\Omega)\big)$ is not dense in 
$C_b(\Omega)$. Let $z\in C(\overline{\mathbb D})$ denote the identity function. Set 
$a=(z,z,\ldots)\in C_b(\Omega)$. Fix $0<\delta<1$, and suppose that 
$f=(f_j)_{j\in\mathbb N}\in C_b(\Omega)$ satisfies $\|f-a\|<\delta$. Then $\|f_j-z\|<\delta$ for every 
$j\in\mathbb N$. We claim that $(1-\delta)\overline{\mathbb D} \subseteq f_j(\overline{\mathbb D})$ for every 
$j\in\mathbb N$. Let $w\in(1-\delta)\overline{\mathbb D}$. For each $j\in\mathbb N$ 
define a function $\phi_j\colon \overline{\mathbb D}\to\mathbb C$ by $\phi_j(\zeta)=w+\zeta-f_j(\zeta)$.
Since
$$|\phi_j(\zeta)| \leq |w|+|\zeta-f_j(\zeta)| < (1-\delta)+\delta= 1\qquad\text{ for every }
\zeta\in\overline{\mathbb D},$$
the function $\phi_j$ maps $\overline{\mathbb D}$ into itself. By Brouwer's fixed-point theorem, there exists 
a number $\zeta_j\in\overline{\mathbb D}$ such that $\phi_j(\zeta_j)=\zeta_j$ yielding 
$w=f_j(\zeta_j)\in f_j(\overline{\mathbb D})$ which proves the claim.

Continuity of $f_j$ implies that each set $f_j(\overline{\mathbb D})$ is connected. Since all these sets contain the 
disk $(1-\delta)\overline{\mathbb D}$, \cite[Theorem~V.1.5]{Dug66} yields that 
$\bigcup_{j\in\mathbb N}f_j(\overline{\mathbb D})$ is connected. Moreover, by \cite[Theorem~V.1.6]{Dug66} its 
closure is also connected, which means that $\sigma(f)$ is connected. Consequently,
$f\notin\mathcal C_2\big(C_b(\Omega)\big)$. Since the open ball with center $a$ and radius $\delta$ does not 
intersect $\mathcal C_2\big(C_b(\Omega)\big)$, the class $\mathcal C_2\big(C_b(\Omega)\big)$ is not dense in 
$C_b(\Omega)$, and hence, none of the spectral classes $\mathcal C_n\big(C_b(\Omega)\big)$ with $n\geq 2$ 
can be dense.

The same proof shows that in the case $\mathbb I=\{1,\ldots,m\}$ none of the sets 
$\mathcal C_n\big(C_b(\Omega)\big)$ with $2\leq n\leq m$ is dense in $C_b(\Omega)$, although these classes 
are nonempty.\qed
\end{example}

\begin{proposition} \label{prop:finite-dimensional-quotient-depth}
Let $\mathscr I$ be a proper closed two-sided ideal in $\mathscr A$ such that the quotient $\mathscr A/\mathscr I$ 
is finite-dimensional, and let 
$$\big(\mathscr A/\mathscr I\big)\big/\rad\big(\mathscr A/\mathscr I\big)\cong
\mathbb M_{n_1}\oplus\cdots\oplus\mathbb M_{n_r}$$
be as in \eqref{eq:wedderburn--artin}. Put $n=n_1+\cdots+n_r$. Then the set 
$\pi^{-1}\big(\mathcal C_n(\mathscr A/\mathscr I)\big)$ is dense in $\mathscr A$ and 
$\pi^{-1}\big(\mathcal C_{n+1}(\mathscr A/\mathscr I)\big)=\varnothing$.
\end{proposition}

\begin{proof}
By Proposition~\ref{prop:depth-finite-dim-algebra} we have $d_\sigma(\mathscr A/\mathscr I)=n$, and therefore 
$\mathcal C_{n+1}(\mathscr A/\mathscr I)=\varnothing$; the second assertion follows.

Denote by $\mathscr D$ the quotient $\big(\mathscr A/\mathscr I\big)\big/\rad\big(\mathscr A/\mathscr I\big)$ 
and let $q\colon\mathscr A/\mathscr I\to\mathscr D$ be the corresponding quotient map. It was shown in the 
proof of Proposition~\ref{prop:depth-finite-dim-algebra} that $\sigma(a)=\sigma(q(a))$ for every 
$a\in\mathscr A/\mathscr I$. Hence, $a\in\mathcal C_n(\mathscr A/\mathscr I)$ if and only if 
$q(a)\in\mathcal C_n(\mathscr D)$, that is, 
$\mathcal C_n(\mathscr A/\mathscr I)=q^{-1}\big(\mathcal C_n(\mathscr D)\big)$.

The set of all matrices in $\mathbb M_m$ with $m$ distinct eigenvalues is dense in $\mathbb M_m$ 
(see, e.g., \cite{Har95}). Applying this to each summand of $\mathscr D$, we conclude that 
$\mathcal C_n(\mathscr D)$ is dense in $\mathscr D$. Since $q\circ\pi\colon\mathscr A\to\mathscr D$ is 
continuous, open and surjective, and 
$\pi^{-1}\big(\mathcal C_n(\mathscr A/\mathscr I)\big)=(q\circ\pi)^{-1}\big(\mathcal C_n(\mathscr D)\big)$, 
Lemma~\ref{lem:open-quotient-maps-density-both-ways} gives the first assertion.
\end{proof}

We now examine what Proposition~\ref{prop:finite-dimensional-quotient-depth} gives for a concrete Banach space 
$X$. Recall that a celebrated construction of Argyros and Haydon \cite{AH11} provides an infinite-dimensional 
separable Banach space $X_{\mathrm{AH}}$ with the \emph{scalar-plus-compact property}. Namely, every operator 
on $X_{\mathrm{AH}}$ is of the form $\lambda I +K$, where $\lambda\in\mathbb C$ and $K$ is compact. 
As usual, let $\mathcal Q(X)$ denote the Calkin algebra $\mathscr B(X)/\mathscr K(X)$ and let 
$\pi\colon\mathscr B(X)\to\mathcal Q(X)$ be the quotient map.

\begin{corollary}
\label{cor:essential-classes-finite-depth}
Let $X=X_{\mathrm{AH}}\oplus\cdots\oplus X_{\mathrm{AH}}$ be the direct sum of $n\in\mathbb N$ copies of the 
Argyros--Haydon space $X_{\mathrm{AH}}$. Then $\mathcal Q(X)\cong\mathbb M_n$, the set
$\pi^{-1}\big(\mathcal C_n(\mathcal Q(X))\big)$ is dense in $\mathscr B(X)$, and 
$\pi^{-1}\big(\mathcal C_{n+1}(\mathcal Q(X))\big)=\varnothing$.
\end{corollary}

\begin{proof}
Every operator $T\in\mathscr B(X)$ is represented uniquely as an $n\times n$ operator matrix 
$T=(T_{rs})_{r,s=1}^{n}$ with $T_{rs}\in\mathscr B(X_{\mathrm{AH}})$. By the scalar-plus-compact property, 
each entry admits a unique decomposition $T_{rs}=\lambda_{rs}I_{X_{\mathrm{AH}}}+K_{rs}$ with 
$\lambda_{rs}\in\mathbb C$ and $K_{rs}\in\mathscr K(X_{\mathrm{AH}})$. Moreover, since $X$ is a finite direct 
sum, the operator $T$ is compact if and only if every entry $T_{rs}$ is compact. Hence, the mapping 
$\Phi\big(\pi(T)\big)=(\lambda_{rs})_{r,s=1}^{n}$ is a well-defined algebra isomorphism from 
$\mathcal Q(X)$ onto $\mathbb M_n$. Since $\mathbb M_n$ is a simple algebra, 
Proposition~\ref{prop:finite-dimensional-quotient-depth} applies with $r=1$ and $n_1=n$.
\end{proof}

The assumption in Proposition~\ref{prop:finite-dimensional-quotient-depth} that the quotient be finite-dimensional 
cannot be omitted, even if the quotient algebra is commutative and semisimple.

\begin{example}\label{ex:nonempty-essential-spectral-class-not-dense}
Let $n\geq2$ and let $\Omega=\{1,\ldots,n\}\times\overline{\mathbb D}$ be the disjoint union of $n$ copies of 
the closed unit disk. By \cite[Theorem~A]{Mot24} there exists an infinite-dimensional complex Banach space $X$ 
such that $\mathcal Q(X)\cong C(\Omega)$ as unital Banach algebras. By 
Proposition~\ref{prop:Ck-spectral-depth} we have $d_\sigma(\mathcal Q(X))=n$, while 
Example~\ref{ex:all-spectral-classes-nonempty-none-dense} shows that $\mathcal C_2(C(\Omega))$ is not dense in 
$C(\Omega)$. Hence, $\mathcal C_2(\mathcal Q(X))$ is not dense in $\mathcal Q(X)$, and 
Proposition~\ref{prop:quotient-Ck-basic-properties} gives that 
$\pi^{-1}\big(\mathcal C_{2}(\mathcal Q(X))\big)$ is not dense in $\mathscr B(X)$.\qed
\end{example}

\section{Perturbation results for operators on Hilbert spaces} \label{sec:perturbations-Hilbert-space}
In Section \ref{sec:topological-properties-spectral-classes}, we proved that for an infinite-dimensional complex 
Banach space $X$, the spectral class $\mathcal C_n(\mathscr B(X))$ is dense in $\mathscr B(X)$. In this section we 
are going to study related perturbation and density results for self-adjoint and normal operators on a 
complex Hilbert space. 

Throughout this section, let $H$ be an infinite-dimensional complex Hilbert space. For a normal operator 
$N\in \mathscr B(H)$ the set $\sigma(N)\setminus\sigma_e(N)$ consists precisely of those isolated points 
$\lambda\in\sigma(N)$ that are eigenvalues of $N$ of finite multiplicity. It follows that a number $\lambda$ belongs 
to $\sigma_e(N)$ if and only if there exists an orthonormal sequence $(e_j)_{j\in\mathbb N}\subseteq H$ such that 
$\lim\limits_{j\to\infty}\|(N-\lambda I)e_j\|=0$. In particular, if $N$ has connected spectrum, then 
$\sigma_e(N)=\sigma(N)$. We refer the reader to \cite[Chapter~XI]{Con90} for details.

The following lemma is our basic perturbation tool.

\begin{lemma}
\label{lem:positive-exact-split}
Let $A\in\mathscr B(H)$ be self-adjoint with connected spectrum and put $b=\max\sigma(A)$. For every $t>0$ 
there exist a positive semi-definite rank-one operator $F\in\mathscr B(H)$ with $\|F\|=t$ and a number $\mu$ 
satisfying $b+\frac{3}{4}t\leq\mu\leq b+t$ such that
$$\sigma(A+F)=\sigma(A)\cup\{\mu\}.$$
Moreover, $\mu$ is a simple eigenvalue of $A+F$ and with respect to the orthogonal decomposition 
$H=H_1\oplus H_1^{\perp}$, where $H_1=\ker(A+F-\mu I)$, we have $A+F=\mu I_{H_1}\oplus A_1$ with 
$\sigma(A_1)=\sigma(A)$.
\end{lemma}

\begin{proof}
Since $\sigma(A)$ is a nonempty compact connected subset of $\mathbb R$, it is either a singleton $\{b\}$ or a 
compact interval $[a,b]$ with $a<b$. 

Suppose first that $\sigma(A)=\{b\}$. Then $A=bI$ by the continuous functional calculus. Choose a unit vector 
$u\in H$ and set $F=t\,(u\otimes u)$. Denote by $[u]$ the linear span of $u$. With respect to the orthogonal 
decomposition $H=[u]\oplus [u]^{\perp}$ we have $A+F=(b+t)I_{[u]}\oplus bI_{[u]^{\perp}}$. It is easily
seen that all the assertions hold with $\mu=b+t$, with $H_1=[u]$, and $A_1=bI_{[u]^{\perp}}$. Note that 
$\sigma(A_1)=\{b\}=\sigma(A)$ because $[u]^{\perp}$ cannot be trivial.

Assume now that $\sigma(A)=[a,b]$ with $a<b$. Let $E_A$ denote the spectral measure of $A$. Since 
$(b-\frac{1}{4}t,b]\cap\sigma(A)\neq\varnothing$, the spectral projection $E_A\big((b-\frac{1}{4}t,b]\big)$ is 
nonzero. Choose a unit vector $u\in E_A\big((b-\frac{1}{4}t,b]\big)H$ and set $F=t\,(u\otimes u)$. Then $F$ is a 
positive semi-definite rank-one operator with $\|F\|=t$. Let $\nu_u$ be the scalar spectral measure given by 
$\nu_u(\Delta)=\langle E_A(\Delta)u,u\rangle$ for Borel sets $\Delta\subseteq\sigma(A)$. Since $\nu_u$ is 
a probability measure supported in $(b-\frac{1}{4}t,b]$, we obtain
$$\langle Au,u\rangle=\int_{(b-\frac{1}{4}t,\,b]}s\,d\nu_u(s)\geq 
\big(b-\tfrac{1}{4}t\big)\nu_u\big((b-\tfrac{1}{4}t,\,b]\big)=b-\tfrac{1}{4}t$$
and therefore
\begin{equation}\label{eq:lower-bound-numerical-range}
\langle (A+F)u,u\rangle=\langle Au,u\rangle+t\geq b+\tfrac{3}{4}t .
\end{equation}
Since $A+F$ is self-adjoint, $\max\sigma(A+F)=\sup\{\langle (A+F)x,x\rangle;\; x\in H,\,\|x\|=1\}$, 
and therefore \eqref{eq:lower-bound-numerical-range} gives $\max\sigma(A+F)\geq b+\frac{3}{4}t$. In particular, 
the spectral projection $Q=E_{A+F}\big((b,\infty)\big)$ of $A+F$ is nonzero, and the subspaces $QH$ and 
$(I-Q)H$ reduce the operator $A+F$.

We claim that the rank of $Q$ is $1$. Suppose, to the contrary, that $\dim(QH)\geq 2$. Then the intersection 
$QH\cap [u]^\perp$ is nontrivial and therefore there is a nonzero $y\in QH$ with $y\perp u$. Denote by 
$\omega_y$ the scalar spectral measure of $A+F$ at $y$. Since $y\in QH$, the measure $\omega_y$ is supported 
in $(b,\infty)$. Hence,
$$\langle (A+F)y,y\rangle-b\|y\|^2=\int_{(b,\infty)}(s-b)\,d\omega_y(s)>0 .$$
On the other hand, $y\perp u$ gives $Fy=0$, and $\sigma(A)\subseteq[a,b]$ gives $A\leq bI$, whence 
$\langle (A+F)y,y\rangle=\langle Ay,y\rangle\leq b\|y\|^2$, a contradiction. Hence, the rank of $Q$ is $1$.
Put $H_1=QH$. Since $\dim(H_1)=1$ and $H_1$ reduces $A+F$, there is a number $\mu\in\mathbb R$ with 
$(A+F)|_{H_1}=\mu I_{H_1}$, and $\mu>b$ by the definition of $Q$. Since $\mu$ is the only spectral value of 
$A+F$ in $(b,\infty)$, we get $\mu=\max\sigma(A+F)\geq b+\frac{3}{4}t$. Moreover, $A\leq bI$ and 
$0\leq F\leq tI$ yield $A+F\leq (b+t)I$ which implies $\mu\leq b+t$. Thus, $b+\frac{3}{4}t\leq \mu\leq b+t$.

It remains to identify $\sigma(A+F)$. Since $F$ has finite rank,
$\sigma_e(A+F)=\sigma_e(A)$, and since $\sigma(A)=[a,b]$ is connected,
the observation at the beginning of this section gives
$\sigma_e(A)=\sigma(A)=[a,b]$. 
Hence, $[a,b]\subseteq\sigma(A+F)$. Furthermore, $A\geq aI$ and $F\geq0$ give $A+F\geq aI$ which implies that 
$\sigma(A+F)\subseteq[a,\infty)$. Since $\sigma(A+F)\cap(b,\infty)=\{\mu\}$, we conclude
$$\sigma(A+F)=[a,b]\cup\{\mu\}=\sigma(A)\cup\{\mu\} .$$
In particular, $\mu$ is an isolated point of $\sigma(A+F)$ with $\ker(A+F-\mu I)=H_1$ one-dimensional, that is, 
$\mu$ is a simple eigenvalue. Let $A_1=(A+F)|_{H_1^{\perp}}$. Then $A+F=\mu I_{H_1}\oplus A_1$ and 
$\sigma(A_1)=\sigma(A+F)\setminus\{\mu\}=\sigma(A)$.
\end{proof}

\begin{theorem}
\label{thm:positive-exact-n}
Let $A\in\mathscr B(H)$ be self-adjoint and suppose that the spectrum $\sigma(A)$ has exactly $k\geq 1$
components. Let $n>k$ be arbitrary. For every $\varepsilon>0$, there exists a positive semi-definite operator 
$F$ of rank $n-k$ such that $\|F\|<\varepsilon$ and $\sigma(A+F)$ has exactly $n$ components.
\end{theorem}

\begin{proof}
Let $\sigma_1,\ldots,\sigma_k$ be the components of $\sigma(A)$ and put $m=n-k$. Let $E$ be the spectral 
measure of $A$. Set $H_j=E(\sigma_j)H$ and $A_j=A|_{H_j}$ for $j=1,\ldots,k$. Then 
$H=H_1\oplus\cdots\oplus H_k$ is an orthogonal decomposition reducing $A$, and $A=A_1\oplus\cdots\oplus A_k$ 
with $\sigma(A_j)=\sigma_j$. Since $H$ is infinite-dimensional, at least one summand is infinite-dimensional; 
we may assume that $\dim(H_k)=\infty$. Denote $\sigma_0=\sigma_k$, $b=\max\sigma_0$, $K_0=H_k$, 
and $G_0=A_k$.

If $k=1$, put $d=\infty$; otherwise let $d=\dist\big(\sigma_0,\sigma(A)\setminus\sigma_0\big)=
\inf\{|s-s'|;\; s\in\sigma_0,\, s'\in\sigma(A)\setminus\sigma_0\}$. Since $\sigma_0$ and 
$\sigma(A)\setminus\sigma_0$ are disjoint compact sets, $d>0$. Choose positive numbers $t_1,\ldots,t_m$ 
with $t_1<\min\{\varepsilon,d\}$ and $t_{j+1}<\tfrac{3}{4}t_j$ for every $1\leq j<m$.

We construct recursively subspaces and operators as follows. Suppose that $G_{j-1}$ is a self-adjoint operator 
on an infinite-dimensional subspace $K_{j-1}\subseteq H_k$ with $\sigma(G_{j-1})=\sigma_0$
(which holds for $j=1$ by the choice of $K_0$ and $G_0$). 
By Lemma~\ref{lem:positive-exact-split} there are a positive semi-definite rank-one operator $F_j$ on $K_{j-1}$ 
with $\|F_j\|=t_j$ and a number $\mu_j$ satisfying 
\begin{equation}\label{eq:mu-j-bounds}
b+\tfrac{3}{4}t_j\leq\mu_j\leq b+t_j.
\end{equation}
The subspace $L_j=\ker(G_{j-1}+F_j-\mu_j I)$ is one-dimensional, which means that 
$K_j=K_{j-1}\cap L_{j}^{\perp}$ is infinite-dimensional, and with respect to the decomposition 
$K_{j-1}=L_j\oplus K_j$ we have 
$$G_{j-1}+F_j=\mu_j I_{L_j}\oplus G_j\qquad\text{where}\qquad  \sigma(G_j)=\sigma_0 .$$
At each step, we extend $F_j$ by zero to $H$. After $m$ steps, we obtain the orthogonal decomposition 
$H_k=L_1\oplus\cdots\oplus L_m\oplus K_m$. Denote $F=F_1+\cdots+F_m$. Then
$$A_k+F=\mu_1I_{L_1}\oplus\cdots\oplus\mu_m I_{L_m}\oplus G_m .$$
The operator $F$ is positive semi-definite. Its rank is $m$, because the subspaces $L_1,\ldots,L_m$ are 
mutually orthogonal and $\ran(F_j)=L_j$. It follows that $\|F\|=\max\{ t_j;\; 1\leq j\leq m\}=t_1<\varepsilon$.
Since $F$ vanishes on $H_1\oplus\cdots\oplus H_{k-1}$, we get 
$$A+F=A_1\oplus\cdots\oplus A_{k-1}\oplus\mu_1I_{L_1}\oplus\cdots\oplus\mu_mI_{L_m}\oplus G_m$$
and therefore $\sigma(A+F)=\sigma_1\cup\cdots\cup\sigma_k\cup\{\mu_1,\ldots,\mu_m\}$. 
It remains to check that this set has exactly $n$ components. The numbers $\mu_1,\ldots,\mu_m$ are pairwise 
distinct. Indeed, by \eqref{eq:mu-j-bounds} and $t_{j+1}<\frac{3}{4}t_j$, we have
$\mu_{j+1}\leq b+t_{j+1}<b+\tfrac{3}{4}t_j\leq\mu_j$. Also, each $\mu_j$ satisfies $\mu_j>b=\max\sigma_0$, 
which gives $\mu_j\notin\sigma_0$. Moreover $\mu_j\leq b+t_j\leq b+t_1<b+d$, whence $\mu_j$ lies at distance 
less than $d$ from $\sigma_0$ and therefore belongs to none of $\sigma_1,\ldots,\sigma_{k-1}$. This shows that 
the components of $\sigma(A+F)$ are precisely the sets $\sigma_1,\ldots,\sigma_k,\{\mu_1\},\ldots,\{\mu_m\}$,
that is, $\sigma(A+F)$ has exactly $k+m=n$ components.
\end{proof}

\begin{corollary}
\label{cor:positive-atmost-n-1}
If $A\in\mathscr B(H)$ is self-adjoint, then for every $\varepsilon>0$ and every $n\in \mathbb N$, there exists a 
positive semi-definite operator $F$ of rank at most $n-1$ such that $\|F\|<\varepsilon$ and $A+F$ belongs to the 
spectral class $\mathcal C_n\big(\mathscr B(H)\big)$.
\end{corollary}

\begin{proof}
If $A\in \mathcal C_n\big(\mathscr B(H)\big)$, then we take $F=0$. Otherwise $\sigma(A)$ has $k<n$
components and Theorem~\ref{thm:positive-exact-n} applies.
\end{proof}

Repeating the recursive construction in the proof of Theorem~\ref{thm:positive-exact-n} indefinitely and choosing 
the norms of the rank-one perturbations to be summable, we obtain a positive semi-definite perturbation belonging 
to the Schatten class $\mathcal S_1(H)$, and therefore to all Schatten classes since 
$\mathcal S_1(H)\subseteq \mathcal S_p(H)$ and $\|T\|_p\leq\| T\|_1$ $(T\in \mathcal S_1(H))$ for every 
$1\leq p\leq \infty$. We use the convention that for $p=\infty$ the class $\mathcal S_p(H)$ equals the ideal of 
compact operators and $\|\cdot\|_\infty$ means the usual operator norm.
 
\begin{theorem}
\label{thm:positive-schatten-infinite-components}
Let $A\in\mathscr B(H)$ be self-adjoint. 
\begin{enumerate}
\item For every $\varepsilon>0$, there exists a positive semi-definite operator $F\in\mathcal S_1(H)$ such that 
$\|F\|_1<\varepsilon$ and $A+F\in\mathcal C_\infty(\mathscr B(H))$.
\item If $\sigma(A)$ has only finitely many components, then $F$ may be chosen so that
$$\sigma(A+F)=\sigma(A)\cup\{\mu_j;\;j\in\mathbb N\},$$
where $(\mu_j)_{j\in\mathbb N}$ is a strictly decreasing sequence of simple isolated eigenvalues converging to 
the right endpoint of a component of $\sigma(A)$.
\end{enumerate}
\end{theorem}

\begin{proof}
If $\sigma(A)$ has infinitely many components, then $A\in\mathcal C_\infty(\mathscr B(H))$ and we may take 
$F=0$. Assume therefore that $\sigma(A)$ has finitely many components $\sigma_1,\ldots,\sigma_k$.
Under this assumption, the construction below proves (i) and (ii) simultaneously.

We use the notation of the proof of Theorem~\ref{thm:positive-exact-n}. Thus, $H=H_1\oplus\cdots\oplus H_k$ is 
the orthogonal decomposition reducing $A$ with $A=A_1\oplus\cdots\oplus A_k$ and $\sigma(A_j)=\sigma_j$, 
where $\dim(H_k)=\infty$, and we write $\sigma_0=\sigma_k$, $b=\max\sigma_0$, $K_0=H_k$ and $G_0=A_k$. 
Let $d=\dist\big(\sigma_0,\sigma(A)\setminus\sigma_0\big)>0$ if $k\geq2$, and $d=\infty$ if $k=1$. Choose 
$t_1>0$ with $\frac{8}{7}t_1<\min\{\varepsilon,d\}$ and put $t_j=8^{1-j}t_1$ for $j\in\mathbb N$, so that
\begin{equation}\label{eq:t-sums}
\sum_{j=1}^{\infty}t_j=\tfrac{8}{7}t_1<\min\{\varepsilon,d\}
\qquad\text{and}\qquad 
t_{j+1}=\tfrac{1}{8}t_j<\tfrac{3}{4}t_j .
\end{equation}

Applying Lemma~\ref{lem:positive-exact-split} recursively as in the proof of 
Theorem~\ref{thm:positive-exact-n}, but without stopping, we obtain positive semi-definite rank-one operators 
$F_j$ with $\|F_j\|=t_j$, one-dimensional subspaces $L_j$, infinite-dimensional subspaces 
$K_0\supseteq K_1\supseteq\cdots$ with $K_{j-1}=L_j\oplus K_j$, self-adjoint operators 
$G_j\in\mathscr B(K_j)$ with $\sigma(G_j)=\sigma_0$, and numbers $\mu_j$ satisfying 
\begin{equation}\label{eq:mu-j-bounds-schatten}
b+\tfrac{3}{4}t_j\leq\mu_j\leq b+t_j ,
\end{equation}
such that $G_{j-1}+F_j=\mu_jI_{L_j}\oplus G_j$. Each $F_j$ is extended by zero to $H$.

Since the operators $F_j$ are positive semi-definite with mutually orthogonal one-dimensional ranges, 
$\|F_j\|_1=t_j$, and \eqref{eq:t-sums} shows that the series $\sum_{j}F_j$ converges in $\mathcal S_1(H)$. 
Let $F$ be its sum. Then $F$ is positive semi-definite and $\|F\|_1=\sum_{j}t_j<\varepsilon$.

Fix $m\in\mathbb N$ and put $R_m=G_m+\sum\limits_{j=m+1}^{\infty}F_j\in\mathscr B(K_m)$. 
Since $F_j$ vanishes on $L_1\oplus\cdots\oplus L_m$ for $j>m$, we have 
\begin{equation}\label{eq:block-form}
(A+F)|_{H_k}=\mu_1I_{L_1}\oplus\cdots\oplus\mu_m I_{L_m}\oplus R_m.
\end{equation}
Moreover $\big\|\sum\limits_{j=m+1}^{\infty}F_j\big\|\leq\sum\limits_{j=m+1}^{\infty}t_j=\frac{1}{7}t_m$, 
and $\sigma(G_m)=\sigma_0$, so
\begin{equation}\label{eq:sigma-Rm}
\sigma(R_m)\subseteq\big\{s\in\mathbb R;\;\dist(s,\sigma_0)\leq\tfrac{1}{7}t_m\big\} .
\end{equation}

We claim that $\sigma\big((A+F)|_{H_k}\big)=\sigma_0\cup\{\mu_j;\;j\in\mathbb N\}$. By 
\eqref{eq:block-form}, each $\mu_j$ is an eigenvalue with eigenvector in $L_j$. Furthermore, $F$ is compact and 
$\sigma(A_k)=\sigma_0$ is connected, so the observations preceding 
Lemma~\ref{lem:positive-exact-split} give $\sigma_e\big((A+F)|_{H_k}\big)=\sigma_e(A_k)=\sigma_0$; 
in particular $\sigma_0\subseteq\sigma\big((A+F)|_{H_k}\big)$. Conversely, let 
$\lambda\notin\sigma_0\cup\{\mu_j;\;j\in\mathbb N\}$ and put $\eta=\dist(\lambda,\sigma_0)>0$. Choose $m$ so 
large that $\frac{1}{7}t_m<\frac{\eta}{2}$. By \eqref{eq:sigma-Rm} every point of $\sigma(R_m)$ is at distance 
at most $\eta/2$ from $\sigma_0$, whence $\lambda\notin\sigma(R_m)$. Since we also have $\lambda\neq\mu_j$ 
$(j=1,\ldots,m)$, the decomposition \eqref{eq:block-form} gives 
$\lambda\notin\sigma\big((A+F)|_{H_k}\big)$, which proves the claim.
As the operator $F$ vanishes on $H_1\oplus\cdots\oplus H_{k-1}$, we conclude
$\sigma(A+F)=\sigma(A)\cup\{\mu_j;\;j\in\mathbb N\}$.

It remains to describe the points $\mu_j$. By \eqref{eq:mu-j-bounds-schatten} and \eqref{eq:t-sums},
$\mu_{j+1}\leq b+t_{j+1}<b+\tfrac{3}{4}t_j\leq\mu_j$.
This means that the sequence $(\mu_j)_{j\in\mathbb N}$ is strictly decreasing, and $t_j\to 0$ gives $\mu_j\to b$. 
Each $\mu_j$ satisfies $\mu_j>b$, hence $\mu_j\notin\sigma_0$. Inequalities $\mu_j\leq b+t_j\leq b+t_1<b+d$ 
imply that $\mu_j$ belongs to no other component of $\sigma(A)$. Consequently, every $\mu_j$ is an isolated point 
of the spectrum $\sigma(A+F)$. To finish the proof, observe that each spectrum $\sigma(R_j)$ lies within distance 
$\frac{1}{7}t_j$ of $\sigma_0$ by \eqref{eq:sigma-Rm}. On the other hand, since $\mu_j\geq b+\frac{3}{4}t_j$, 
we have $\mu_j\notin\sigma(R_j)$. The equality \eqref{eq:block-form} with $m=j$ shows that 
$\ker(A+F-\mu_jI)=L_j$. Thus, each $\mu_j$ is a simple isolated eigenvalue of $A+F$. This proves that 
$\sigma(A+F)$ has infinitely many components, that is, $A+F\in\mathcal C_\infty(\mathscr B(H))$.
\end{proof}

Since the perturbation $F$ in Theorem~\ref{thm:positive-schatten-infinite-components} is of trace class, it is 
compact, and therefore $\sigma_e(A+F)=\sigma_e(A)$. Thus, an arbitrarily small positive perturbation can turn a 
self-adjoint operator with connected spectrum into an operator whose spectrum has infinitely many components, 
while leaving the essential spectrum unchanged. The points $\mu_j$ produced by the construction lie in 
$\sigma(A+F)\setminus\sigma_e(A+F)$. These points are precisely the isolated eigenvalues of finite multiplicity 
that distinguish $\sigma(A+F)$ from $\sigma(A)$. This illustrates the invariance recorded in 
Proposition~\ref{prop:quotient-Ck-basic-properties}, applied to the ideal $\mathscr K(H)$ in $\mathscr B(H)$.
\medskip

Let $\mathscr B(H)_{\mathrm{sa}}$ denote the real Banach space of all self-adjoint operators on $H$. 
Theorem~\ref{thm:Cn-density} shows that $\mathcal C_n(\mathscr B(X))$ is dense in $\mathscr B(X)$
for every infinite-dimensional complex Banach space $X$, but this gives no information about approximation within 
$\mathscr B(H)_{\mathrm{sa}}$. The results of this section fill this gap.

\begin{corollary}
\label{cor:selfadjoint-residual}
For every $n\in\mathbb N$, the set $\mathcal C_n\big(\mathscr B(H)\big)\cap\mathscr B(H)_{\mathrm{sa}}$ is 
open and dense in $\mathscr B(H)_{\mathrm{sa}}$, and consequently 
$\mathcal C_\infty\big(\mathscr B(H)\big)\cap\mathscr B(H)_{\mathrm{sa}}$ is a dense $G_\delta$-subset of 
$\mathscr B(H)_{\mathrm{sa}}$. The same holds within the set of all positive semi-definite operators on $H$.
\end{corollary}

\begin{proof}
Openness is immediate from Theorem~\ref{thm:Ck-open}, since the intersection of an open subset of 
$\mathscr B(H)$ with $\mathscr B(H)_{\mathrm{sa}}$ is relatively open. Density follows from 
Corollary~\ref{cor:positive-atmost-n-1}: given $A\in\mathscr B(H)_{\mathrm{sa}}$ and $\varepsilon>0$, the 
operator $A+F$ produced there is self-adjoint, lies in $\mathcal C_n(\mathscr B(H))$, and satisfies 
$\|F\|<\varepsilon$. As $\mathscr B(H)_{\mathrm{sa}}$ is a Baire space, the intersection
$$\mathcal C_\infty\big(\mathscr B(H)\big)\cap\mathscr B(H)_{\mathrm{sa}}
=\bigcap_{n=1}^{\infty}\Big(\mathcal C_n\big(\mathscr B(H)\big)\cap\mathscr B(H)_{\mathrm{sa}}\Big)$$
is a dense $G_\delta$-subset. The perturbation $F$ of Corollary~\ref{cor:positive-atmost-n-1} is positive 
semi-definite, which means that the operator $A+F$ is positive semi-definite whenever $A$ is. Since the positive
semi-definite operators form a closed subset of $\mathscr B(H)_{\mathrm{sa}}$, this set is a complete metric space.
Thus, we may apply the same arguments as before, and the last assertion holds.
\end{proof}

Our next goal is to prove a perturbation theorem for normal operators. To this end, we will need the following 
basic topological lemma.

\begin{lemma}
\label{lem:boundary-infinite-compact-set}
If $\Delta\subseteq\mathbb C$ is a compact infinite set, then $\partial \Delta$ contains a point which is 
an accumulation point of $\partial \Delta$ and therefore of $\Delta$.
\end{lemma}

\begin{proof}
Suppose that $\partial\Delta$ is finite. Then
$\mathbb C\setminus\partial\Delta$ is connected, and
\[
\mathbb C\setminus\partial\Delta
=
\Delta^{\circ}\cup(\mathbb C\setminus\Delta)
\]
is a disjoint union of open sets. Hence, one of these sets must be empty.
If $\Delta^{\circ}=\varnothing$, then $\Delta=\partial\Delta$, so $\Delta$
is finite, contrary to the assumption. If $\mathbb C\setminus\Delta=\varnothing$,
then $\Delta=\mathbb C$, contradicting the compactness of $\Delta$. Thus,
$\partial\Delta$ is infinite. Being an infinite compact subset of $\mathbb C$, it has an
accumulation point. Since $\partial\Delta$ is closed, this accumulation point
belongs to $\partial\Delta$, and, as $\partial\Delta\subseteq\Delta$, it is
also an accumulation point of $\Delta$.
\end{proof}

We are now ready to prove the normal-operator perturbation theorem.

\begin{theorem}\label{thm:normal-compact-perturbation-Cinfty}
Let $N\in\mathscr B(H)$ be normal. Then, for every $\varepsilon>0$, there exists
a compact operator $K\in\mathscr B(H)$ such that
$\|K\|<\varepsilon$, the operator $N+K$ is normal, and
$N+K\in\mathcal C_\infty(\mathscr B(H))$.
\end{theorem}

\begin{proof}
Let $E$ be the spectral measure of $N$ and, for $x\in H$, let $\nu_x$ be the scalar spectral measure given by 
$\nu_x(\Delta)=\langle E(\Delta)x,x\rangle$. We first construct a point $\lambda_0\in\partial\sigma(N)$ and
a separable infinite-dimensional reducing subspace $H_0$ of $N$ such that $\lambda_0\in\sigma_e(N_0)$, 
where $N_0=N|_{H_0}$.

Suppose first that $\sigma(N)$ is finite. Then $\sigma(N)=\partial\sigma(N)$ and $E$ is purely atomic. One has an
orthogonal decomposition $H=\bigoplus\limits_{\mu\in\sigma(N)}\ran\big(E(\{\mu\})\big)$. Since $H$ is 
infinite-dimensional, there is $\lambda_0\in\sigma(N)$ with $\dim\big(\ran\big(E(\{\lambda_0\})\big)\big)=\infty$. 
Let $H_0$ be any separable infinite-dimensional subspace of $\ran\big(E(\{\lambda_0\})\big)$. As an eigenspace of 
$N$, the subspace $\ran\big(E(\{\lambda_0\})\big)$ reduces $N$. The same holds for $H_0$; moreover 
$N_0=\lambda_0I_{H_0}$, whence $\sigma_e(N_0)=\{\lambda_0\}$.

If $\sigma(N)$ is infinite, then we can apply Lemma~\ref{lem:boundary-infinite-compact-set}. Hence, there is a 
point $\lambda_0\in\partial\sigma(N)$ which is an accumulation point of $\sigma(N)$. Choose distinct points 
$\alpha_j\in\sigma(N)$ with $\alpha_j\to\lambda_0$ and $\alpha_j\ne\lambda_0$. Since the $\alpha_j$ are 
distinct, we may choose radii $r_j>0$ so small that the open disks $U_j=B(\alpha_j,r_j)$ are pairwise disjoint 
and $r_j\leq|\alpha_j-\lambda_0|$. Since $\alpha_j\in U_j\cap\sigma(N)$, the projection $E(U_j)$ is nonzero.
Let $u_j\in\ran\big(E(U_j)\big)$ be a unit vector. The disjointness of the sets $U_j$ implies that 
$(u_j)_{j\in\mathbb N}$ is orthonormal.

Put $H_j=\overline{\{f(N)u_j;\; f\in C(\sigma(N))\}}$ and $H_0=\bigoplus_{j=1}^{\infty}H_j$. Each $H_j$ is 
invariant under $f(N)$ for every $f\in C(\sigma(N))$. In particular, it is invariant under $N$ and $N^*$, which means
that $H_j$ and $H_0$ reduce $N$. Since $\sigma(N)$ is a compact metric space, the algebra $C(\sigma(N))$ 
is separable. Every mapping $f\mapsto f(N)u_j$ is continuous and therefore each $H_j$ is separable. It follows that
$H_0$ is separable. Taking $f\equiv 1$ shows $u_j\in H_j\subseteq H_0$ and therefore $N_0u_j=Nu_j$.
For $z\in U_j$ we have $|z-\lambda_0|\leq|z-\alpha_j|+|\alpha_j-\lambda_0|<2|\alpha_j-\lambda_0|$. Since 
$\nu_{u_j}$ is supported in $U_j$, we have
$$\|(N_0-\lambda_0I)u_j\|^2=\int_{U_j}|z-\lambda_0|^2\,d\nu_{u_j}(z)\leq 4|\alpha_j-\lambda_0|^2,$$
which gives $\|(N_0-\lambda_0I)u_j\|\to 0$. As $(u_j)_{j\in\mathbb N}$ is orthonormal, the characterization of the 
essential spectrum recalled at the beginning of this section gives $\lambda_0\in\sigma_e(N_0)$.

By the Weyl--von Neumann--Berg theorem \cite[Theorem~39.4]{Con00} there is a diagonalizable operator 
$D\in\mathscr B(H_0)$ with $D-N_0$ compact and $\|D-N_0\|<\varepsilon/2$. Let $(v_m)_{m\in\mathbb N}$ 
be an orthonormal basis of $H_0$ with $Dv_m=\gamma_mv_m$. Since $D-N_0$ is compact and $H_0$ reduces $N$,
$\sigma_e(D)=\sigma_e(N_0)\subseteq\sigma(N_0)\subseteq\sigma(N)$; in particular, $\lambda_0\in\sigma_e(D)$. 
The spectral measure of $D$ is 
$E_D(\Delta)=\sum\limits_{m:\,\gamma_m\in\Delta}\langle\,\cdot\,,v_m\rangle v_m$. It follows that 
$\dim\big(E_D\big(B(\lambda_0,r)\big)H_0\big)=\infty$ for every $r>0$. Hence, the set 
$\{m\in\mathbb N;\;\gamma_m\in B(\lambda_0,1/j)\}$ is infinite for every $j\in\mathbb N$, and we may choose 
recursively a strictly increasing sequence $(m_j)_{j\in\mathbb N}$ with $\gamma_{m_j}\in B(\lambda_0,1/j)$. 
In particular $\gamma_{m_j}\to\lambda_0$.

The set $\sigma(D)\setminus\sigma(N)$ is at most countable. Indeed, if $\mu\in\sigma(D)\setminus\sigma(N)$, 
then $\mu\notin\sigma_e(D)$. This means that $\mu$ is an isolated point of $\sigma(D)$. Since $\sigma(D)$ is 
the closure of $\{\gamma_m;\;m\in\mathbb N\}$, every isolated point of $\sigma(D)$ lies in this countable set.

Fix $j\in\mathbb N$ and put $\tau_j=\min\{\varepsilon\,2^{-j-2},\,1/j\}$. Since $\lambda_0\in\partial\sigma(N)$, 
the set $B(\lambda_0,\tau_j)\setminus\sigma(N)$ is nonempty and open, hence uncountable. Hence, it is 
not contained in the countable set $\sigma(D)\setminus\sigma(N)$. We may thus choose 
$\mu_j\in B(\lambda_0,\tau_j)\setminus\big(\sigma(N)\cup\sigma(D)\big),$
and, shrinking $\tau_j$ further if necessary, we may assume that the $\mu_j$ are pairwise distinct. Then 
$|\mu_j-\lambda_0|<\varepsilon\,2^{-j-2}$ and
$|\mu_j-\gamma_{m_j}|\leq|\mu_j-\lambda_0|+|\lambda_0-\gamma_{m_j}|\longrightarrow 0$
since $|\gamma_{m_j}-\lambda_0|<1/j$. Discarding finitely many indices, we may assume 
$|\mu_j-\gamma_{m_j}|<\varepsilon/2$ for every $j\in\mathbb N$.

Define a diagonal operator $R\in\mathscr B(H_0)$ by 
$$Rv_{m_j}=(\mu_j-\gamma_{m_j})v_{m_j}\quad \text{if } j\in\mathbb N\qquad\text{and}\qquad Rv_m=0 
\quad\text{if } m\notin\{m_j;\;j\in\mathbb N\} .$$
Since $\mu_j-\gamma_{m_j}\to0$, the operator $R$ is compact, and its norm is $\|R\|\leq \varepsilon/2$. 
Put $K=\big((D-N_0)+R\big)\oplus 0$ with respect to the decomposition $H=H_0\oplus H_0^{\perp}$. 
Then $K$ is compact and 
$\|K\|\leq\|D-N_0\|+\|R\|<\varepsilon$. Writing $N=N_0\oplus N_1$ we obtain $N+K=(D+R)\oplus N_1$,
which is normal, being the direct sum of a diagonalizable operator and the restriction of $N$ to a reducing subspace.

Since $R$ is compact, $\sigma_e(D+R)=\sigma_e(D)\subseteq\sigma(N)$. Fix $j\in\mathbb N$. By construction 
$(D+R)v_{m_j}=\mu_jv_{m_j}$, i.e., $\mu_j$ is an eigenvalue of the normal operator $D+R$. Since 
$\mu_j\notin\sigma(N)$, this eigenvalue is not in $\sigma_e(D+R)$. An eigenvalue of a normal operator lying 
outside its essential spectrum is an isolated point of the spectrum. Hence, $\mu_j$ is isolated in $\sigma(D+R)$. 
Since $\sigma(N)$ is closed and $\mu_j\notin\sigma(N)$, there exists $\delta_j>0$ such that the punctured disk 
$\{z;\;0<|z-\mu_j|<\delta_j\}$ meets neither $\sigma(D+R)$ nor $\sigma(N)$. Since 
$\sigma(N_1)\subseteq\sigma(N)$ and $\sigma(N+K)=\sigma(D+R)\cup\sigma(N_1)$,
this punctured disk misses $\sigma(N+K)$. Therefore, $\mu_j$ is an isolated point of $\sigma(N+K)$, that is, 
$\{\mu_j\}$ is a component of $\sigma(N+K)$. As the points $\mu_j$ are pairwise distinct, $\sigma(N+K)$ has 
infinitely many components. This proves that $N+K\in\mathcal C_\infty(\mathscr B(H))$.
\end{proof}

The operator $K$ in Theorem~\ref{thm:normal-compact-perturbation-Cinfty} is in general not normal. Indeed, the 
following proposition shows that for a unitary operator without eigenvalues the only compact normal perturbation 
which preserves normality is the trivial one.

\begin{proposition}\label{prop:unitary-normal-perturbation-trivial}
Let $U\in\mathscr B(H)$ be a unitary operator and suppose that it has no eigenvalues. 
If $K$ is a compact normal operator such that $U+K$ is normal, then $K=0$.
\end{proposition}

\begin{proof}
Set $A=U^*K$. The normality of $U+K$ gives $(U+K)^*(U+K)=(U+K)(U+K)^*$
which simplifies to $U^*K+K^*U=UK^*+KU^*$ as $K$ and $U$ are normal. Since 
$$ U(A+A^*)U^*=UAU^*+UA^*U^*=KU^*+UK^*=U^*K+K^*U=A+A^*,$$
the unitary operator $U$ commutes with the compact self-adjoint operator $A+A^*$. If $A+A^*\neq 0$, then it has 
a nonzero eigenvalue whose eigenspace is finite-dimensional and invariant under $U$, which is a contradiction. 
Hence, $A^*=-A$. Since $K=UA$, normality of $K$ yields $AA^*=A^*A=UAA^*U^*$, which means that $U$ 
commutes with the compact self-adjoint operator $AA^*$. By the same argument applied to $AA^*$, we get 
$AA^*=0$. It follows from $\|A\|^2=\|AA^*\|=0$ that $A=0$ and consequently $K=UA=0$. 
\end{proof} 

Let $\mathbb T$ be the unit circle and let $m$ be the normalized Lebesgue measure on $\mathbb T$.
Consider the multiplication operator $M_z$ given by $(M_zf)(z)=zf(z)$ on the complex Hilbert space 
$L^2(\mathbb T,m)$. It is well-known that $M_z$ is unitary and has no eigenvalues.

\begin{corollary}\label{cor:no-normal-nonzero-perturbation}
There is no nonzero compact normal operator $K$ such that $M_z+K$ is normal. 
\end{corollary}

\begin{proof}
Since $M_z$ is a unitary operator without eigenvalues, Proposition~\ref{prop:unitary-normal-perturbation-trivial}
applies.
\end{proof}

Let $\mathcal Q(H)$ be the Calkin algebra of $H$ and $\pi$ the corresponding quotient map.

\begin{proposition}\label{prop:essential-components-normal-operator}
Let $N\in\mathscr B(H)$ be a normal operator and let $n\in\mathbb N$. If 
$N\in\pi^{-1}\big(\mathcal C_n(\mathcal Q(H))\big)$, then $N\in\mathcal C_n(\mathscr B(H))$.
\end{proposition}

\begin{proof}
At the beginning of Section~\ref{sec:perturbations-Hilbert-space}, we already observed that the set 
$\sigma(N)\setminus\sigma_e(N)$ consists of isolated eigenvalues of $N$ of finite multiplicity. Hence, every 
component of $\sigma_e(N)$ is also a component of $\sigma(N)$.
\end{proof}

The normality assumption cannot be omitted in the previous proposition, because for a general operator 
$T\in \mathscr B(H)$ the set $\sigma(T)\setminus\sigma_e(T)$ need not consist of isolated points.

\begin{example}\label{ex:essential-components-do-not-imply-spectral-components}
Let $H_n$ be the orthogonal sum of $n\geq2$ copies of $\ell^2$. Let $S$ be the unilateral shift on $\ell^2$, and 
choose numbers $0<r_1<r_2<\cdots<r_n$. Define the operator $T\in\mathscr B(H_n)$ by
$T=r_1S\oplus r_2S\oplus\cdots\oplus r_nS$.
It is well known that $\sigma(S)=\overline{\mathbb D}$ and $\sigma_e(S)=\mathbb T$, whence 
$\sigma(r_jS)=r_j\overline{\mathbb D}$ and $\sigma_e(r_jS)=r_j\mathbb T$ for every $j=1,\ldots,n$. 
Since the spectrum of a direct sum is the union of the spectra of its summands and $r_1<\cdots<r_n$, we obtain 
$\sigma(T)=r_n\overline{\mathbb D}$. This shows that $T\notin\mathcal C_2(\mathscr B(H_n))$. On the other hand, 
it is well-known that a direct sum of operators is Fredholm if and only if each summand is Fredholm. Hence,
the essential spectrum of $T$ is $\sigma_e(T)=\bigcup_{j=1}^{n}r_j\mathbb T$, the union of $n$ pairwise disjoint 
circles. This shows that $T\in \pi^{-1}\big(\mathcal C_n(\mathcal Q(H_n))\big)$.\qed
\end{example}


\section{Rank-one perturbations of quasinilpotent compact shifts} \label{sec:quasinilpotent-rank-one-perturbation}
Let $X$ be an infinite-dimensional complex Banach space. In this section we consider a class of
compact quasinilpotent forward shifts on $X$ and study their rank-one perturbations. Our starting
point is a result of Tcaciuc \cite[Proposition~2.4]{Tca19}, which describes the point spectrum of
rank-one perturbations of arbitrary quasinilpotent operators. We first specialize this result to the
compact shifts considered here and express its conclusion in terms of the spectral class
$\mathcal C_\infty(\mathscr B(X))$. We then show that, for a suitable choice of the rank-one
perturbation, the exceptional scalar allowed by Tcaciuc's general theorem can be eliminated
completely.
 
Let $(e_j)_{j\in\mathbb N}$ be any basic sequence in $X$ (as in Corollary~\ref{cor:spectral-depth-BX}, we refer 
the reader to \cite[Section~1.5]{AK06} for the details). Without loss of generality, we may assume that 
the sequence is normalized. Let $Y=\overline{\Span}\{e_j;\;j\in\mathbb N\}$, and let 
$(\psi_j)_{j\in\mathbb N}\subseteq Y^*$ be the corresponding sequence of biorthogonal functionals. For each 
$j\in\mathbb N$, let $\xi_j\in X^*$ be a Hahn--Banach extension of $\psi_j$, i.e., $\|\xi_j\|=\|\psi_j\|$. Moreover,
$\|\xi_j\|\leq 2K$ for every $j\in\mathbb N$, where $K$ is the basis constant of $(e_j)_{j\in\mathbb N}$.
Throughout this section, the above sequences $(e_j)_{j\in\mathbb N}\subseteq X$ and 
$(\xi_j)_{j\in\mathbb N}\subseteq X^*$ are fixed.

Let $w=(w_j)_{j\in\mathbb N}\in\ell^1$ be a sequence of positive numbers. Then the series
$$\sum_{j=1}^\infty w_j(e_{j+1}\otimes \xi_j)$$
converges absolutely in $\mathscr B(X)$, since 
$\sum\limits_{j=1}^{\infty} w_j\|e_{j+1}\otimes \xi_j\|=\sum\limits_{j=1}^{\infty} w_j\|\xi_j\|\leq 2K\|w\|_1$, 
and therefore defines an operator $S_w$, called a \emph{weighted forward shift}. Since the 
sequence of finite-rank partial sums of this series converges in norm to $S_w$, the operator is compact.

\begin{lemma} \label{lem:weighted-compact-quasinilpotent-shift}
Let $w=(w_j)_{j\in\mathbb N}\in\ell^1$ be a nonincreasing sequence of positive numbers. Then $S_w$ is 
quasinilpotent if and only if
\begin{equation} \label{eq:weighted-compact-quasinilpotent-shift}
\lim_{k\to\infty}\big(w_1\cdots w_k\big)^{1/k}=0.
\end{equation}
\end{lemma}

\begin{proof}
Let $W_0=1$ and for every $k\in\mathbb N$ let $W_k=w_1\cdots w_k$. Since 
$\langle e_{j+1},\xi_i\rangle=\delta_{i,j+1}$, a simple induction gives
$$S_{w}^{\,k}x=\sum_{j=1}^{\infty}\Big(\prod_{i=j}^{j+k-1}w_i\Big)\langle x,\xi_j\rangle\,e_{j+k}\qquad 
(x\in X,\ k\in\mathbb N).$$
In particular $S_{w}^{\,k}e_1=W_k\,e_{k+1}$ which implies $\|S_{w}^{\,k}\|\geq W_k$
as the sequence $(e_j)_{j\in\mathbb N}$ is normalized. If $S_w$ is quasinilpotent, then 
$W_{k}^{1/k}\leq\|S_{w}^{\,k}\|^{1/k}\to r(S_w)=0$, which proves 
\eqref{eq:weighted-compact-quasinilpotent-shift}. Conversely, assume 
\eqref{eq:weighted-compact-quasinilpotent-shift}. Since $w$ is nonincreasing, we have 
$w_{j+t}\leq w_t$ for all $j,t\in\mathbb N$, whence
$$\prod_{i=j}^{j+k-1}w_i=w_j\prod_{t=1}^{k-1}w_{j+t}\leq w_j\,W_{k-1}\qquad (j, k\in\mathbb N).$$
Together with $\|\xi_j\|\leq 2K$ this yields
$$\|S_w^{\,k}\|\leq\sum_{j=1}^{\infty}\Big(\prod_{i=j}^{j+k-1}w_i\Big)\|\xi_j\|
\leq 2K\,W_{k-1}\sum_{j=1}^{\infty}w_j=2K\|w\|_1\,W_{k-1} .$$
By \eqref{eq:weighted-compact-quasinilpotent-shift}, $W_{k-1}^{1/k}=\big(W_{k-1}^{1/(k-1)}\big)^{(k-1)/k}\to0$
and therefore
\[r(S_w)=\lim_{k\to\infty}\|S_{w}^{\,k}\|^{1/k}
\leq\lim_{k\to\infty}\big(2K\|w\|_1\big)^{1/k}W_{k-1}^{1/k}=0 . \qedhere\]
\end{proof}

From now on, let $w$ be a nonincreasing sequence of positive numbers in $\ell^1$ satisfying 
\eqref{eq:weighted-compact-quasinilpotent-shift}, and denote $W_k=w_1\cdots w_k$ for 
$k\in\mathbb N$ and $W_0=1$. For $x\in X$ and $\varphi\in X^*$ we define a formal power series
$$G_{x,\varphi}(z)=\sum_{k=0}^{\infty}\langle S_{w}^{\,k}x,\varphi\rangle\, z^{k+1}.$$
Since $|\langle S_{w}^{\,k}x,\varphi\rangle|^{1/(k+1)}\leq
\big(\|S_{w}^{\,k}\|\,\|x\|\,\|\varphi\|\big)^{1/(k+1)}\to 0$, because $S_w$ is quasinilpotent, the radius of 
convergence of this series is infinite by the Cauchy--Hadamard formula. Hence, $G_{x,\varphi}$ is an entire 
function. 

The connection between the function $G_{x,\varphi}$ and rank-one perturbations is the rank-one resolvent criterion
used by Tcaciuc in \cite[Lemma~2.1]{Tca19}. In our notation it takes the form of the following lemma. For the sake of completeness, we include its proof.

\begin{lemma}\label{lem:G_n-entire-function}
Let $x\in X$ and $\varphi\in X^*$. A nonzero number $\lambda$ belongs to the spectrum 
$\sigma(S_w+x\otimes\varphi)$ if and only if $G_{x,\varphi}(1/\lambda)=1$.
\end{lemma}

\begin{proof}
Since the spectrum of the operator $I-y\otimes \xi$ is the set $\{1,1-\langle y,\xi\rangle\}$, we see that 
$I-y\otimes \xi$ is invertible if and only if $\langle y,\xi\rangle\neq 1$. Pick 
$\lambda\in\mathbb C\setminus\{0\}$. Since $S_w$ is quasinilpotent, by Lemma 
\ref{lem:weighted-compact-quasinilpotent-shift}, the operator $\lambda I-S_w$ is invertible. Hence, the identity
$$\lambda I-(S_w+x\otimes\varphi)=(\lambda I-S_w)\big(I-(\lambda I-S_w)^{-1}x\otimes\varphi\big)$$
implies that $\lambda\in\sigma(S_w+x\otimes\varphi)$ if and only if the operator 
$I-(\lambda I-S_w)^{-1}x\otimes\varphi$ is not invertible, that is, if and only if
$\langle(\lambda I-S_w)^{-1}x,\varphi\rangle=1$. Since the spectral radius of $S_w$ is zero, expansion of 
$(\lambda I-S_w)^{-1}$ into its Neumann series gives 
$$(\lambda I-S_w)^{-1}x=\sum_{k=0}^{\infty}\lambda^{-k-1}S_{w}^{\,k}x .$$
It follows that $\lambda\in\sigma(S_w+x\otimes\varphi)$ if and only if
\[G_{x,\varphi}(1/\lambda)=\sum_{k=0}^{\infty}\langle S_{w}^{\,k}x,\varphi\rangle\,\lambda^{-k-1}
=\langle(\lambda I-S_w)^{-1}x,\varphi\rangle=1 .\qedhere\]
\end{proof}

We now apply the general rank-one perturbation theorem of Tcaciuc. Recall that
\cite[Proposition~2.4]{Tca19} and its proof show that if $T$ is quasinilpotent and
$F=x\otimes\varphi$ is rank one such that
$$
\langle T^kx,\varphi\rangle\neq0
$$
for infinitely many $k$, then for every nonzero scalar $\alpha$, with at most one exception,
the point spectrum of $T+\alpha F$ is countably infinite. In our setting compactness
turns this conclusion immediately into a statement about spectral components.

\begin{theorem}\label{thm:Cinfty-never-finite-rank-stable}
Let $w\in\ell^1$ be a nonincreasing sequence of positive numbers satisfying 
\eqref{eq:weighted-compact-quasinilpotent-shift}. Let $x\in X$ and $\varphi\in X^*$ be such that 
$\langle S_{w}^{\,k}x,\varphi\rangle\neq 0$ for infinitely many $k\geq 0$. Then, for the rank-one operator 
$F=x\otimes\varphi$, we have
$$S_w+\alpha F\in\mathcal C_\infty\big(\mathscr B(X)\big)$$
for every $\alpha\in\mathbb C\setminus\{0\}$ with at most one exception.
\end{theorem}

\begin{proof}
We apply \cite[Proposition~2.4]{Tca19} to $T=S_w$ and
$F=x\otimes\varphi$. By assumption,
$$
\langle T^kx,\varphi\rangle
=
\langle S_w^{\,k}x,\varphi\rangle
\neq0
$$
for infinitely many $k$. Hence, for every nonzero $\alpha$, with at most one exception,
$\sigma_p(S_w+\alpha F)$ is countably infinite. Since $S_w+\alpha F$ is compact, every nonzero point of its 
spectrum is isolated. Hence, each singleton $\{\lambda_j\}$ is a component of $\sigma(S_w+\alpha F)$. 
Therefore, $\sigma(S_w+\alpha F)$ has infinitely many components, that is, 
$S_w+\alpha F\in\mathcal C_\infty(\mathscr B(X))$.
\end{proof}

The hypothesis of Theorem~\ref{thm:Cinfty-never-finite-rank-stable} is easily satisfied. Indeed, take $x=e_1$. 
Then $S_{w}^{\,k}e_1=W_k e_{k+1}$, which gives the entire function
\begin{equation}
\label{eq:G-e1-phi}
G_{e_1,\varphi}(z)=\sum_{k=1}^{\infty}W_{k-1}\langle e_k,\varphi\rangle\, z^k
\end{equation}
and the condition of the theorem then reads $\langle e_k,\varphi\rangle\neq 0$ for infinitely 
many $k\in\mathbb N$, because $W_{k-1}>0$ for every $k$. Let $(\gamma_j)_{j\in\mathbb N}\in \ell^1$ be
such that $\gamma_j\ne 0$ for infinitely many $j$. Since $\|\xi_k\|\leq 2K$, the series 
$\sum_{k=1}^{\infty}\gamma_k \xi_k$ converges absolutely in $X^*$; its sum $\varphi$ satisfies 
$\langle e_k,\varphi\rangle=\gamma_k$ for every $k\in\mathbb N$, and therefore 
$\langle e_k,\varphi\rangle\ne 0$ for infinitely many $k$.

The possible exceptional scalar in Theorem~\ref{thm:Cinfty-never-finite-rank-stable} cannot in general
be omitted. Indeed, Tcaciuc \cite[Example]{Tca19} gives a compact quasinilpotent weighted shift $T$ and a rank-one operator $F$ such that $T-F$ is quasinilpotent, whereas $\sigma_p(T+\alpha F)$ is countably infinite for every $\alpha\in\mathbb C\setminus\{0,-1\}$. Thus, the exceptional nonzero scalar allowed by \cite[Proposition~2.4]{Tca19} may actually
occur.

We now show that, for every shift $S_w$ considered here, the
rank-one perturbation $F$ can be chosen in such a way that no nonzero
exceptional scalar occurs. For this purpose, let
\[
\Gamma
=
\left\{
(\gamma_k)_{k\in\mathbb N}\in\ell^1:
\gamma_1=0,\;
\gamma_{2k}=0,\;
\gamma_{2k+1}\neq 0
\text{ for every } k\in\mathbb N
\right\}.
\]
For each $\gamma=(\gamma_k)_{k\in\mathbb N}\in\Gamma$, define
\[
\varphi_\gamma
=
\sum_{k=1}^{\infty}\gamma_k\xi_k
\qquad\text{and}\qquad
F_\gamma
=
e_1\otimes\varphi_\gamma.
\]
Since $(\gamma_k)_{k\in\mathbb N}\in\ell^1$, the series defining
$\varphi_\gamma$ converges absolutely in $X^*$, so
$\varphi_\gamma\in X^*$. Thus, $F_\gamma$ is a bounded rank-one
operator. The next result
shows that every such choice gives a perturbation for which all nonzero
scalars are nonexceptional.

\begin{theorem}
\label{thm:Cinfty-gamma-sequence-criterion}
Let $w\in\ell^1$ be a nonincreasing sequence of positive numbers
satisfying \eqref{eq:weighted-compact-quasinilpotent-shift}. Then, for
every $\gamma\in\Gamma$, the operator $F_\gamma$ is nilpotent and
\[
S_w+\alpha F_\gamma
\in
\mathcal C_\infty\bigl(\mathscr B(X)\bigr)
\]
for every $\alpha\in\mathbb C\setminus\{0\}$.
\end{theorem}

\begin{proof}
Let $\gamma=(\gamma_k)_{k\in\mathbb N}\in\Gamma$. Since
$\langle e_1,\varphi_\gamma\rangle=\gamma_1=0$, we have
\[
F_\gamma^2
=
\langle e_1,\varphi_\gamma\rangle F_\gamma
=
0,
\]
so $F_\gamma$ is nilpotent.

By \eqref{eq:G-e1-phi}, we have
\[
G_{e_1,\varphi_\gamma}(z)
=
\sum_{k=1}^{\infty}
W_{2k}\gamma_{2k+1}z^{2k+1}.
\]
Since $W_{2k}>0$ and $\gamma_{2k+1}\neq 0$ for every
$k\in\mathbb N$, this entire function has infinitely many nonzero
coefficients. By the Great Picard theorem, there is at most one number
$\beta\in\mathbb C$ for which the equation
$G_{e_1,\varphi_\gamma}(z)=\beta$ does not have infinitely many
solutions. We claim that such an exceptional value, if it exists, must
be $0$. Indeed, suppose that $\beta\neq 0$ were exceptional. Since
$G_{e_1,\varphi_\gamma}$ is odd, the equation
$G_{e_1,\varphi_\gamma}(z)=-\beta$ has infinitely many solutions if and
only if the equation $G_{e_1,\varphi_\gamma}(z)=\beta$ has infinitely
many solutions. Thus, $-\beta$ would also be exceptional. Since
$\beta\neq-\beta$, this would give two exceptional values, contradicting
the Great Picard theorem. Therefore, the only possible exceptional value
is $0$.

Let now $\alpha\in\mathbb C\setminus\{0\}$. Since $1/\alpha\neq 0$,
the equation
\[
G_{e_1,\varphi_\gamma}(z)
=
\frac{1}{\alpha}
\]
has infinitely many solutions. Let $(\mu_j)_{j\in\mathbb N}$ be a
sequence of distinct such solutions. Since the zeros of the nonzero
entire function $G_{e_1,\varphi_\gamma}-1/\alpha$ have no accumulation
point in $\mathbb C$, the sequence $(\mu_j)_{j\in\mathbb N}$ is
unbounded. Passing to a subsequence, we may assume that
$|\mu_j|\to\infty$. Since
$G_{e_1,\varphi_\gamma}(0)=0\neq 1/\alpha$, we also have
$\mu_j\neq 0$ for every $j\in\mathbb N$.

Set $\lambda_j=1/\mu_j$. Then the numbers $\lambda_j$ are distinct,
nonzero, and belong to $\sigma(S_w+\alpha F_\gamma)$, and
$\lambda_j\to 0$. Since $S_w+\alpha F_\gamma$ is compact, every
nonzero point of its spectrum is isolated. Hence each singleton
$\{\lambda_j\}$ is a component of
$\sigma(S_w+\alpha F_\gamma)$. Therefore,
$\sigma(S_w+\alpha F_\gamma)$ has infinitely many components, that is,
\[
S_w+\alpha F_\gamma
\in
\mathcal C_\infty\big(\mathscr B(X)\big).\qedhere
\]
\end{proof}

Since $c\gamma\in\Gamma$ for every $\gamma\in\Gamma$ and every
$c\in\mathbb C\setminus\{0\}$, we have
$\varphi_{c\gamma}=c\varphi_\gamma$ and
$F_{c\gamma}=cF_\gamma$. Consequently, the perturbations obtained by
Theorem~\ref{thm:Cinfty-gamma-sequence-criterion} may be chosen with
arbitrarily small norm.

\begin{corollary}
\label{cor:affine-rank-one-line-Cinfty}
Let $X$ be an infinite-dimensional complex Banach space. Then there exist a compact quasinilpotent
operator $S\in\mathscr K(X)$ and a rank-one nilpotent operator $F\in\mathscr B(X)$ such that
$$
\{S+\alpha F;\;\alpha\in\mathbb C\setminus\{0\}\}
\subseteq
\mathcal C_\infty\big(\mathscr B(X)\big).
$$
Moreover, $F$ may be chosen with arbitrarily small norm.
\end{corollary}

\subsection*{Acknowledgements}
The first author was supported by the Slovenian Research and Innovation Agency program P2-0268.
The second author was supported by the Slovenian Research and Innovation Agency program P1-0222 and grant J1-50002.

\subsection*{Data Availability Statement}

No datasets were generated or analysed during the current study.

\subsection*{Conflicts of Interest}

The authors declare that they have no conflicts of interest relevant to the content of this article.

\subsection*{AI Declaration}
The authors used ChatGPT 5.6 Sol during the development of this work to assist with identifying relevant references, examining 
mathematical arguments, and language editing. All mathematical arguments were independently checked and substantially revised by 
the authors, who take full responsibility for the content and the final version of the manuscript.


\end{document}